\documentclass[a4paper,12pt]{article} 
\usepackage[latin1]{inputenc}
\usepackage{pdfsync}
\usepackage{amsmath}
\usepackage{amssymb}
\newtheorem{defi}{Definition}[section]
      
      \newtheorem{prop}[defi]{Proposition}
      \newtheorem{prop-def}[defi]{Proposition-Definition}
      \newtheorem{def-prop}[defi]{Definition-Proposition}
      \newtheorem{thm}[defi]{Theorem}
      
      \newtheorem{lem}[defi]{Lemma}
      \newtheorem{cor}[defi]{Corollary}

\newcommand {\br}[1]{\bar #1}
\newcommand{\cqfd}{\hfill $\Box$}
\newcommand {\esp}[1]{\underset{ #1}{\mathbb E}}
\newcommand {\OO}[1]{\mathcal{O}/\pi^{ #1}\mathcal{O}}

\def\O{\mathcal{O}}
\def\unif{\pi}

\def\eps{\varepsilon}

\def\R{{\mathbb R}}
\def\Z{{\mathbb Z}}
\def\N{{\mathbb N}}
\def\C{{\mathbb C}}
\def\CC{{\mathcal C}}

\def\E{\mathcal{E}}

\def\E{\mathcal{E}}

\def\pro{\mathrm{p}}
\def\F{\mathbb{F}}

\def\L{\mathcal{L}}

\def\gg{\mathfrak {g}}

\def\z{\mathcal{O}}

\def\p{\pi}

\newcommand{\s}[1]{\langle #1 \rangle}

\RequirePackage{graphics}

\begin{document}

\title{Strong Banach property (T)\footnote{This text is an improved version of the
author's master thesis (mémoire), under the supervision of Vincent Lafforgue,
at Université Paris 7, Denis Diderot.}\\ \large after Vincent Lafforgue}
\author{Benben Liao \footnote{Institut de mathématiques de
Jussieu, 7C16, 175 rue du Chevaleret, 75013 Paris, liaob@math.jussieu.fr,
liao@hotmail.co.uk.}}
\maketitle

{\bf Title in French} Propriété (T) renforcée banachique.

{\bf Abstract in French} Ce texte reprend, avec plus de d\'etails et en simplifiant une preuve
dans la section 3, les parties de \cite{duke} et \cite{JTA} traitant des groupes
p-adiques. On y prouve que $SL_3$ sur un corps local non archimédien $ F $ a la
propriété  (T) renforcée banachique. Les applications sont les suivantes : tout
groupe algébrique connexe presque $F$-simple $G$ sur $F$ dont l'algèbre de
Lie contient $sl_{3}(F)$ a la propriété (T) renforcée banachique,   les suites
d'expanseurs construites à partir d'un réseau de $ G (F) $ ne se plongent
uniformément dans aucun espace de Banach de type $>1$, et toute action affine
isométrique de $ G (F) $ sur un espace de Banach de type $> 1 $ admet un point
fixe.     



\tableofcontents

\addcontentsline{toc}{section}{Introduction}

\vskip 10mm

\noindent {\bf Introduction}

\vskip 3mm
This text takes, with more details and simplifying a proof in section 3, the
parts of \cite{duke} and \cite{JTA} treating p-adic groups. 

A locally compact group $G$ is said to have Kazhdan's property (T) if there exists an idempotent $\pro \in C^{*}_{max}(G)$
, such that for any unitary representation $(H,\pi)$ of $G$, the image of
$\pi(\pro)$ is exactly the space of $G$-invariant vectors $H^{G}$. In the
articles \cite{duke} and \cite{JTA} of Vincent Lafforgue, a strong property
(T) is proposed in the sense that the space of representation is allowed to be a
Banach space of non trivial type, and the representation is not
necessarily isometric but has a small exponential growth.                     

In this text, Banach spaces are always complex Banach spaces with very precise norms, which are not equivalent in general.

\begin{defi}\label{defi-classe-type}
We say that a class of Banach spaces $\E$ is of type $>1$ if the following
equivalent conditions are satisfied:
\begin{itemize}
\item there exists $n \in \N^{*}$ and $\epsilon >0$ such that any Banach space $E$ in $\E$ does not contain $\ell^{n}_{1}$ $(1+\epsilon)$-isometrically,
\item there exists $p>1$ (called the type) and $T \in \R_{+}$ such that for any $E$ in $\E$, $n \in \N^{*}$ and $x_{1},...,x_{n} \in E$, we have
\begin{eqnarray}\label{ineg-type}\Big(\esp{\epsilon_{i}=\pm 1} \|\sum _{i=1}^{n}\epsilon_{i}x_{i}\|_{E}^{2}\Big)^{\frac{1}{2}}\leq T\Big( \sum_{i=1}^{n} \|x_{i}\|^{p}_{E}\Big)^{\frac{1}{p}}.\end{eqnarray}
\end{itemize}
\end{defi}

In this definition we should say "uniformly of type $>1$" instead of "of type
$>1$". We refer to \cite{maurey} for a review about the notion of type. We
simply recall that $\ell^{n}_{1}$ means $\C^{n}$ endowed with the $\ell_{1}$ norm and
that a Banach space contains $\ell^{n}_{1}$ $(1+\epsilon)$-isometrically if
there exists $i: \ell_1^{n} \to E$ such that $\|x\| \leq \|i(x)\| \leq
(1+\epsilon)\|x\|$ for any $x \in \ell^{n}_{1}$. We say that a class $\E$ is
stable under duality if for any $E$ in $\E$, $E^{*}$ is in $\E$. We say that a
class $\E$ is stable under complex conjugation if for any $E$ in $\E$, $\bar{E}$
is in $\E$. Any class of type $>1$ is included in a class of type $>1$ which is
stable under duality and complex conjugation.                      

Let $G$ be a locally compact group, and $dg$ a left invariant Haar measure,
which endows $C_{c}(G)$ with a structure of $\C$-algebra by convolution. We call
length over $G$ a continuous function $\ell:G \to \R_{+}$ satisfying
$\ell(g^{-1})=\ell(g)$  and $\ell(g_{1}g_{2}) \leq \ell(g_{1}) + \ell(g_{2})$
for $g_{1},g_{2} \in G$.    

If $\ell$ is a length over $G$, and $\E$ is a class of Banach spaces, we note
$\E_{G,\ell}$ the class of continuous representations $(E,\pi)$ of $G$ in the
Banach space $E\in\E$ such that $\|\pi(g)\|\leq e^{\ell(g)}$ for
any $g \in G$. By $\CC^{\E}_{\ell}(G)$ we denote the completion of
$C_{c}(G)$ with respect to the norm $\|f\|=\sup_{(E,\pi)\in
\E_{G,\ell}}\|\pi(f)\|_{\L(E)}.$ If $\E$ is stable under duality and complex
conjugation, $\CC^{\E}_{\ell}(G)$ is a Banach algebra endowed with complex
conjugation $\int_{G} f(g)e_{g}dg \mapsto \int_{G} \bar f(g)e_{g}dg $ and the
usual involution $f\mapsto f^*$ defined by $\int_{G} f(g)e_{g}dg \mapsto
\int_{G}  e_{g^{-1}}\bar f(g)dg $, since for any $(E,\pi)\in \E_{G,\ell}$, the
conjugate representation $(\bar E, g\mapsto \bar \pi(g))$, the contragredient
representation $(E^{*},g\mapsto {}^{t}\pi(g^{-1}))$, and the conjugate
contragredient representation $(\bar E^{*},g\mapsto {}^{t}\bar \pi(g^{-1}))$
also belong to $\E_{G,\ell}$.                

\begin{defi}\label{renf-ban}
Let $G$ be a locally compact group. We say that $G$ has strong Banach property
(T) if for any class $\E$ of type $>1$, stable under duality and complex
conjugation, and for any length $\ell$ over $G$, there exists $s_0>0$, such that
for any $C \in \R_{+}$ and any $0< s\leq s_0$, there exists a real and
self-adjoint idempotent $\pro$ in $\CC^{\E}_{C+s\ell}(G)$ such that for any representation $(E,\pi) \in
\E_{G,C+s\ell}$, the image of $\pi(\pro)$ is the subspace of $E$ consisting of
all $G$-invariant vectors.                 
\end{defi}

\noindent {\bf Remark.} In the conditions of the definition, such $\pro$ is
invariant under $$\int_{G} f(g)e_{g}dg \mapsto \int_{G} e_{g^{-1}}f(g)dg $$
(which is the composition of the complex conjugation and the usual involution),
therefore for $(E,\pi) \in \E_{G,C+s\ell}$, $\mathrm{Ker} \pi(\pro)$ is the
orthogonal complement of $\mathrm{Im} {}^{t}\pi(p)=(E^{*})^{G}$. One deduces immediately
that such $\pro$ is unique and that it is a central element of
$\CC^{\E}_{C+s\ell}(G)$.          

\noindent {\bf Remark.} Let $\ell$ be a length over $G$ and $\E$ a class of
Banach spaces such that $\CC^{\E}_{\ell}(G)$ has an idempotent $\pro$ such that
for any representation $(E,\pi) \in \E_{G,\ell}$, the image of $\pi(\pro)$ is
the subspace of $E$ consisting of all $G$-invariant vectors. Suppose that $G$ is
not compact and that the class $\E$ is defined by a super-property (see the
paragraph 2 of \cite{maurey} for this notion). Then the class $\E$ is of type
$>1$. In fact, for any surjective morphism $E \to F$ in the category
$\E_{G,\ell}$ (F being endowed with the quotient norm of that of $E$), $E^{G}
\to F^{G}$ must be surjective since $E^{G}$ and $F^{G}$ are the images of the
idempotent $\pro$ as the subspaces of $E$ and $F$. But the morphism $L^{1}(G)
\to \C,f \mapsto \int_{G} f dg$ is surjective and $G$ is not compact, $L^{1}(G)$
does not have nonzero $G$-invariant vector (precisely $G$ acts on $L^{1}(G)$ by
the left regular representation). Therefore $L^{1}(G)$ does not belong to the
class $\E$. As the class $\E$ is defined by a super-property, it is of type
$>1$. Therefore, whenever the group $G$ is not compact, one cannot expect a
statement stronger than the definition above, if the class $\E$ is defined by a
super-property.                              

The purpose of this text is to prove the following result.

\begin{thm}\label{main}
Let $G$ be a connected almost $F$-simple algebraic group over a non archimedian
local field $F$. If its Lie algebra $\gg$ contains $sl_{3}(F)$, $G$ has strong
Banach property (T).               
\end{thm}

We prove as in \cite{duke} that cocompact lattice inherits strong Banach
property (T). For non cocompact lattice, it is still true if we limit ourselves
to isometric representations, as in the following definition.               

\begin{defi}
Let $G$ be a locally compact group. We say that $G$ has Banach property (T) if for any class $\E$ of type $>1$, stable under duality and complex conjugation, there exits a real and self-adjoint idempotent $\pro$ in $\CC^{\E}_{0}(G)$ such that for any representation $(E,\pi) \in \E_{G,0}$, the image of $\pi(\pro)$ consists of exactly $G$-invariant vectors in $E$.
\end{defi}

Strong Banach property (T) evidently implies Banach property (T). As in \cite{duke} we show that Banach property (T) is inherited by any lattice.

We deduce the following statement.

\begin{cor}
Let $F$ be a non archimedian local field and $G$ a connected, almost $F$-simple algebraic group over $F$ whose Lie algebra contains $sl_{3}(F)$. Let $\Gamma$ be a lattice of $G(F)$ and $\Gamma_{n}$ a series of subgroups of finite indices of $\Gamma$ such that the indices tends to infinity. We fix a finite system of generators of $\Gamma$. Then the Cayley graph of $\Gamma/\Gamma_{n}$ is a series of expanders which cannot be imbedded uniformly into any Banach space of type $>1$, more precisely for any Banach space $E$ of type $>1$ there exists $C\in \R_{+}$ such that for any $n \in \N$ and for any $1$-Lipschitz function $f:\Gamma/\Gamma_{n} \to E$ of zero average, $\esp{x \in \Gamma/\Gamma_{n}} \|f(x)\|_{E}^{2}\leq C$.
\end{cor}

During a discussion with Uri Bader and inspired by an idea of Bader and Gelander formulated in the remark 3.3 of \cite{gelander-note}, Vincent Lafforgue realized that strong Banach property (T) implies Banach property (F) in the following sense.

\begin{defi}\label{def-f}
Let $G$ be a locally compact group. We say that $G$ has Banach property (F) if any continuous isometric affine action of $G$ on a Banach space of type $>1$ admits a fixed point.
\end{defi}

We deduce the following corollary, which strengthens remark 1.7(2) of \cite{uri}
and some particular cases of theorem B of \cite{uri} and of theorem 3.1 of
\cite{gelander-note}.       

\begin{cor}\label{prop-f}
Let $F$ be a non archimedian local field and $G$ a connected, almost $F$-simple algebraic group over $F$ whose Lie algebra contains $sl_{3}(F)$. Then $G$ has Banach property (F).
\end{cor}

Here is a list of open problems from \cite{JTA}.

\noindent a)
Show that for any non archimedian local field $F$, $Sp_{4}(F)$ has strong Banach property (T). This implies that any almost simple algebraic group over a non archimedian local field, of split rank $\geq 2$, has strong Banach property (T).

\noindent b)
Prove strong Banach property (T) for almost simple algebraic group over an
archimedian local field, of split rank $\geq 2$, in particular for $SL_{3}(\R)$.
In fact in \cite{duke}, a hilbertian variant of strong property (T) for
$SL_{3}(\R)$ is proved. With the notations of paragraph 2 of \cite{duke}, strong
Banach property (T) is implied by the following statement: for any class $\E$ of
type $>1$, there exist $\alpha \in \R^{*}_{+}$ and $C \in \R_{+}$ such that for
any $\epsilon \in (-1,1)$ and any $E \in \E$, we have
$$\| (T_{0}-T_{\epsilon})\otimes 1\|_{\L(L^{2}(S^{2},E))}\leq C|\epsilon|^{\alpha}. $$

\noindent c)
Show that the family of Ramanujan graphs obtained by quotient of an arithmetic
cocompact subgroup of $SL_{2}(\R)$ by the arithmetic subgroups of finite indices
does not admit a uniform imbedding in Banach spaces of type $>1$. This question
is difficult since $SL_{2}(\R)$ does not have property (T) and the reason for
these graphs are expanders is that $SL_{2}(\R)$ has property ($\tau$), i.e. the
trivial representation is isolated among the $L^{2}(SL_{2}(\R)/\Gamma)$ with
$\Gamma$ being an arithmetic group. However, we do not know if for $\E$ a class
of type $>1$, the trivial representation of $SL_{2}(\R)$ is isolated among
$L^{2}(SL_{2}(\R)/\Gamma,E)$ with $\Gamma$ an arithmetic group and $E\in \E$. In
fact the proof of property ($\tau$) uses arithmetics. The most elementary proof
of the fact that these graphs are expanders is the one explained in
\cite{davidoff} but it is based on the calculation of traces and we do not know
how to adapt it to coefficients in Banach spaces of type $>1$ (or even uniformly
convex).                            

\noindent d)
Show that any family of expanders considered above does not admit a uniform
imbedding into Banach spaces of finite cotype (see paragraph 3 of \cite{pisier}
for a discussion on this problem).               

In his article,  Lafforgue says he was inspired by discussions with Uri Bader,
Assaf Naor and Gilles Pisier.   

{\bf Acknowledgement.}
I would like to thank Vincent Lafforgue, for his patience in answering my
questions concerning his articles and correcting my text, and also his
recommendation of books and tri-semesters at IHP on related subjects. I also
thank George Skandalis and Stéphane Vassout for their kind help and advice on
administrative matters during the year of my master study, which means a lot to
a foreign student in France like me. Finally, I would like to thank my tutor, Xiaonan Ma for his
advice and recommendation of operator algebra team of IMJ, and Paris Graduate
School of Mathematical Sciences, for providing the master scholarship to me to
study in Paris.  


\section{Review of fast Fourier transform}
Here we follow paragraph 1 in $\cite{JTA}$. We recall the general setting of
fast Fourier transform, following some notes of Jean-Fran\c cois Mestre.       

 Let $G$ be a finite abelian group. We recall that the Fourier transform $T_{G}:\ell^{2}(G) \rightarrow \ell^{2}(\hat{G})$ is given by $T_{G}f(\chi)=\esp{x \in G} \chi(x)f(x)$. This normalization, where $T_{G}$ is not an isometry, is not customary but more convenient for us in the following discussion.

 Let $H$ be a subgroup of $G$. The principle of fast Fourier transform
 is to write $T_{G}$ as the product of two matrices with blocks and then the
 blocks are essentially the matrices of $T_{H}$ and $T_{G/H}$. The idea is to
 decompose the average over $x \in G$ into an average over the classes modulo
 $H,$ and an average over $G/H$. We choose a section $\sigma : G/H \rightarrow
 G$ of the projection $\pi: G \rightarrow G/H$, i.e. $\pi \circ \sigma =
 \mathrm{Id}_{G/H}$.
 Therefore $T_{G}=T_{2,G,H} \circ T_{1,G,H}$, where $T_{1,G,H}:\ell^{2}(G)
 \rightarrow \ell^{2}((G/H) \times \hat{H})$ is given by $$f\mapsto \Big[
 (x',\chi ')\mapsto \esp {x,\pi(x)=x'  } \chi'(x-\sigma(x'))f(x)\Big],$$ and
 $T_{2,G,H}:\ell^{2}((G/H)\times \hat H) \to \ell^{2}(\hat G)$ is given by
 $$g\mapsto  \Big[\chi\mapsto \esp{x'\in G/H} \chi(\sigma(x'))g(x',\chi_{\big |
 H})\Big].$$                

We see that $T_{1,G,H}$ is block diagonal: the blocks are indexed by $G/H$ and
equal to $T_{H}$. Similarly $T_{2,G,H}$ is block diagonal: the blocks are
indexed by $\hat{H}$ and obtained from $T_{G/H}$ by multiplying the columns of
this matrix by a complex number of module 1: more precisely the choice of a
section $\sigma':\hat{H} \rightarrow \hat{G}$ of the restriction $\hat{G}
\rightarrow \hat{H}$ identifies the block of $T_{2,G,H}$ indexed by $\chi' \in
\hat{H}$ with the matrix obtained from $T_{G/H}$ by multiplying the column
indexed by $x' \in G/H$ by the unit $\sigma'(\chi')(\sigma(x'))$.                   

Let $A$ be a finite matrix, where rows and columns are indexed by finite set
$I$ and $J$, and let $E$ be a Banach space. We have operator $A\otimes
1_{E}:\ell^{2}(J,E)\to \ell^{2}(I,E)$, where $\ell^{2}(I,E)$ is the
space of functions on $I$ with values in $E$ with the norm
$\|f\|=\sqrt{\esp{x\in I}\|f(x)\|_{E}^{2}}$ and so is $\ell^2(J,E)$. It is clear
that the norm of $A \otimes 1_{E}$ does not change if we multiply $A$ by a unit,
or we multiply the rows or the columns of $A$ by some units. 

We have then for any Banach space $E$,
 $$\|T_{1,G,H}\otimes 1_{E}\|
 _{\L(\ell^{2}(G,E), \ell^{2}((G/H)\times \hat H,E))}
      =\|T_{H}\otimes 1_{E}\|_{\L(\ell^{2}(H,E),\ell^{2}(\hat H,E))}$$ $$\text{\ \  and\ \  }
      \|T_{2,G,H}\otimes 1_{E}\|_{\L( \ell^{2}((G/H)\times \hat H,E),\ell^{2}(\hat G,E))}
      =\|T_{G/H}\otimes 1_{E}\|_{\L(\ell^{2}(G/H,E),\ell^{2}(\widehat {G/H},E))}.$$
We deduce the following proposition, where we replace the chain of length 2 of
finite abelian groups $0\subset H\subset G$ with a chain of length $n:$
$0=G_{0}\subset G_{1}\subset \dots \subset G_{n}=G$.            

      \begin{prop}\label{fast} for any Banach space $E$,
      $$\|T_{G}\otimes1_{E}\|_{\L(\ell^{2}(G,E),\ell^{2}(\hat G,E))}\leq \prod_{i=1}^{n} \|T_{(G_{i}/G_{i-1})}\otimes1_{E}\|_{\L(\ell^{2}((G_{i}/G_{i-1}),E),\ell^{2}(\widehat {(G_{i}/G_{i-1})},E))}.$$
      \end{prop}

The interest of fast Fourier transform in computer science is to reduce the time
of computations. Take $G=\Z/2^{n}\Z$ as an example. The number of operations for
a computation of the image under $T_{G}$ of a function in
$\ell^{2}(G)$ is $2^{2n}$, but the decomposition of $T_{G}$ into the product of
$n$ matrices with blocks (of size $2 \times 2$), which is associated with the
chain $0=G_{0}\subset G_{1}\subset \dots \subset G_{n}=G$ where
$G_{i}=2^{n-i}\Z/2^{n}\Z$, allows a computation with $O(n2^{n})$ operations.

\section{Fourier transform and type}
We follow paragraph 2 in \cite{JTA}. We recall the relationship between type and the operation $T_{G}$ considered in the previous paragraph. First it is evident that for any finite abelian group G and for any Banach space $E$ we have $\|T_{G}\otimes 1_{E}\|\leq 1$, since in the matrix of $T_{G}$ the sum of the absolute values of coefficients over a row or a column is equal to 1 (the regular norm of $T_{G}$, in the sense of \cite{pisier}, is therefore equal to $1$). First $\|T_{G}\otimes 1_{\ell^{1}(G)}\|= 1$ because the vector $(x,y)\mapsto \delta_{x,y}$ of $\ell^{2}(G,\ell^{1}(G))$ is of norm $1$ (since the norm $\ell^{2}$ is defined by an average and the norm $\ell^{1}$ by a sum) and its image is also of norm $1$. Therefore if $\E$ is a class of Banach spaces defined by a super-property one cannot expect to have $\sup_{E\in \E}\|T_{G}\otimes 1\|<1$ unless the class $\E$ is of type $>1$. The following proposition is theorem 1 of \cite{bourgain-pacific}(note remark 1 at the end of \cite{bourgain-pacific}).

 \begin{prop}
 Let  $\E$ be a class of type $>1$. There exists $p>1$ and $M\in \R_{+}$ such that for any compact abelian group G, any Banach space $E$ in the class $\E$ and any sequence with finite support $(x_{\gamma})_{\gamma \in \hat G}$ of elements in $E$, we have (with $dg$ the Haar measure over $G$ of mass $1$)
 $$ \Big( \int_{g\in G} \big\|\sum _{\gamma\in \hat G}x_{\gamma} \gamma(g)\big\|_{E}^{2}dg\Big)^{1/2} \leq M\Big(\sum _{\gamma\in \hat G}\|x_{\gamma}\|^{p}\Big)^{1/p}.
 $$
 \end{prop}

 We only use this statement for finite group $G$. We can suppose $p\in (1,2]$ in the statement. By replacing the sum with average, the inequality here becomes
   $$ \Big( \int_{g\in G} \big\|\esp{\gamma\in \hat G}x_{\gamma} \gamma(g)\big\|_{E}^{2}dg\Big)^{1/2} \leq M
  (\sharp G)^{\frac{1}{p}-1}
  \Big(\esp{\gamma\in \hat G}\|x_{\gamma}\|^{p}\Big)^{1/p}.
 $$
The inequality of H\" older shows that for all family $(x_{\gamma})_{\gamma \in \hat G}$ of elements in $E$ we have
   $$\Big(\esp{\gamma\in \hat G}\|x_{\gamma}\|^{p}\Big)^{1/p} \leq
  \big(\esp{\gamma\in \hat G}\|x_{\gamma}\|^{2}\big)^{\frac{1}{2}}.$$
By permuting the role of $G$ and $\hat G$, we then obtain the following corollary.

  \begin{cor}\label{bourgain-holder}
  Let $\E$ be a class of type $>1$. There exists $\epsilon >0$ and $C \in \R_{+}$ such that for any finite abelian group $G$ and any Banach space $E$ in the class $\E$, we have
 $$\|T_{G}\otimes 1_{E}\|\leq C(\sharp G)^{-\epsilon}.$$
  \end{cor}

  Let $F$ be a non archimedian local field, $\O$ its ring of integers and $\pi$ a uniformizer of $\O$.

  The following two corollaries are immediate consequences of corollary \ref{bourgain-holder}.

  \begin{cor} \label{OF-alpha} For any class $\E$ of type $>1$, there exists $h \in \N^{*}$ and $\alpha >0$ such that for any
  $E\in \E$, $\|T_{\O/\pi^{h}\O}\otimes 1_{E}\|\leq e^{-\alpha}$.
  \end{cor}

  \begin{cor} \label{OF-beta} For any class $\E$ of type $>1$, there exists $C>0$ and $\beta >0$ such that for any
  $E\in \E$ and $n \in \N^{*}$, $\|T_{\O/\pi^{n}\O}\otimes 1_{E}\|\leq Ce^{-\beta n}$.
  \end{cor}

  In the rest of this paragraph we make precise the logic relation between the statements of corollaries \ref{OF-alpha} and \ref{OF-beta}. It is clear that the statement of corollary \ref{OF-beta} implies that of corollary \ref{OF-alpha}. Thanks to the fast Fourier transform, and more precisely thanks to proposition \ref{fast}, the statement of corollary \ref{OF-alpha} implies that of corollary \ref{OF-beta}. In fact, suppose the statement of corollary \ref{OF-alpha} holds for $h$. Let $n=ah+b$ with $b \in \{0,...,h-1\}$ and $a \in \N$. We consider the chain of subgroups
  $$0\subset   \pi^{ah}\O/\pi^{n}\O\subset   \pi^{(a-1)h}\O/\pi^{n}\O \subset  \dots \subset   \pi^{h}\O/\pi^{n}\O    \subset \O/\pi^{n}\O, $$
  we see thanks to proposition \ref{fast} that $\|  T_{\O/\pi^{n}\O}\otimes 1_{E}\|\leq \Big(\|  T_{\O/\pi^{h}\O}\otimes 1_{E}\|\Big)^{a}$. The statement of corollary \ref{OF-beta} follows, by taking $\beta=\frac{\alpha}{h}$ and $C=e^{\alpha}$.

\section{A variant of fast Fourier transform}
This paragraph is an improvement of the proof of lemma 3.2 in paragraph 3 of \cite{JTA}.

Let $G$, $\br{G}$ be two abelian groups, $Z$ a finite set, and $\kappa$ a bilinear map from $G \times \br{G}$ to the set of unitary operators on $\ell^2(Z)$ whose image consists of commuting matrices, i.e. $$\kappa(x + x',y)= \kappa(x,y) \kappa(x',y),$$ $$ \kappa(x,y+ y')=\kappa(x,y) \kappa(x,y'),$$ $$\kappa(x,0)=1_{\ell^2(Z)},\kappa(0,y)=1_{\ell^2(Z)},$$ and  $$\kappa(x,y)\kappa(x',y')\kappa(x,y)^{-1}\kappa(x',y')^{-1}=1,$$ for any $\ x, x' \in G, y, y' \in \br{G}$. Define $T_{G,\br{G},Z}$ as follows:
\begin{eqnarray} \label{def-T}
  T_{G,\br{G},Z}: \ell^{2}(G\times Z) &\rightarrow & \ell^{2}(\br{G}\times
  Z)\nonumber\\ 
   f &\mapsto & \Big[(y,z) \mapsto \esp{x \in G} \big(\kappa
   (x,y)f(x,\cdot)\big)(z)\Big],
\end{eqnarray}

where $\kappa(x,y)f(x,\cdot)$ is the image by $\kappa(x,y)$ (acting on
$\ell(Z)$) of $f(x,\cdot).$

Let G and $\br{G}$ contain the following subgroups $$0=G_{0} \subset G_{1} \subset G_{2} ... \subset G_{n-1} \subset G_{n}=G$$
$$0=\br{G}^{n} \subset \br{G}^{n-1} \subset \br{G}^{n-2} ... \subset \br{G}^{1}
\subset \br{G}^{0}=\bar G,$$ such that for any $i=0,...,n,$ we have
$\kappa(x,y)=1_{\ell^2(Z)}$, for $x \in G_{i}, y \in \br{G}^{i}.$ Define for any $i=0,...,n-1$, the map
$\kappa_i$ from $ G_{i+1}/G_{i} \times \br{G}^{i}/ \br{G}^{i+1}$ to the set of
unitary operators on $\ell^2(Z)$, by $\kappa_i([x_{i}],[y_{i}])=
\kappa(x_{i},y_{i})$, for any $ x_{i} \in G_{i+1},y_{i} \in \br{G}^{i}$, where
$[\cdot]$ are the corresponding projections $G_{i+1} \rightarrow G_{i+1}/G_{i},
\br{G}_{i} \rightarrow \br{G}_{i}/\br{G}_{i+1}$. Define $T_{G_{i+1}/G_{i},
\br{G}^{i}/ \br{G}^{i+1},Z}$ as in (\ref{def-T}) with $\kappa$ replaced by
$\kappa_i$.



For any finite set $S$, any operator $A\in\mathcal{L}(\ell^2(S))$ and any Banach
space $E$, we have defined the operator norm $\|A\otimes 1_E\|_{\mathcal
L(\ell^2(Z,E))}$ in paragraph 1.
With these settings, we have the following variant of fast Fourier transform.    

\begin{prop}\label{VFFT}
With $Z,\kappa,T_{G,\br{G},Z}$ and the two chains of abelian subgroups
satisfying the conditions above, for any Banach space $E$ satisfying
$\|\kappa(x,y)\otimes 1_E\|_{\mathcal L(\ell^2(Z,E))}\leq 1$ for any $x\in G,y\in \bar G$, we have the
following.
\begin{itemize}
\item (i) $$\|T_{G,\br{G},Z} \otimes 1_{E}\| \leq \prod_{i=0}^{n-1} \| T_{G_{i+1}/G_{i}, \br{G}^{i}/ \br{G}^{i+1},Z}\otimes 1_{E}\|.$$
\item (ii) Let $I \subset \{0,1,..,n-1\}$ be a subset. If for any $i\in I$,
$\kappa_i$ is scalar valued, i.e.
    $\kappa_i(x,y)=\lambda_i (x,y)\mathrm{Id}_{\ell^2(Z)},$ 
     and $\lambda_i$ is a non degenerate (i.e. if $\lambda_i(x,y)=1$ for any
     $x\in G_{i+1}$ then $y=0,$ and if $\lambda_i(x,y)=1$ for any
     $y\in \bar G^{i}$ then $x=0$ ) bilinear form on $G_{i+1}/G_i \times
     \bar G^i /\bar G^{i+1}$ with values of complex numbers of norm $1,$
    then $$\|T_{G,\br{G},Z} \otimes 1_{E}\|_{\L(\ell^2(G\times
    Z,E),\ell^2(\bar G\times Z,E))} $$$$\leq \prod_{i \in I}
    \|T_{G_{i+1}/G_{i}}\otimes
    1_{E}\|_{\L\big(\ell^2(G_{i+1}/G_i,E),\ell^2(\widehat{G_{i+1}/G_i},E)\big)},$$
    where $T_{G_{i+1}/G_i}:\ell^2(G_{i+1}/G_i)\to\ell^2(\widehat{G_{i+1}/G_i})$ is the
    usual Fourier transform defined in paragraph 1.
\end{itemize} 
\end{prop}


We prove (i) by induction. 
Suppose $n=2.$ We use notations $0 \subset H=G_{1} \subset G, 0 \subset
\br{H}=\br{G}_{1} \subset \br{G},$ and we want to prove $$\| T_{G,\br{G},Z}
\otimes 1_{E}\| \leq \| T_{H,\br{G}/\br{H},Z} \otimes 1_{E}\| \cdot \| T_{G/H,\br{H},Z} \otimes 1_{E}\|.$$

First for each of the two projections
$$ p: G \rightarrow G/H$$
$$\pi: \br{G} \rightarrow \br{G}/\br{H},$$

we choose a section respectively
$$s: G/H \rightarrow G$$
$$\sigma: \br{G}/\br{H} \rightarrow \br{G}.$$

We then define the following two operators,
$$\begin{matrix} T_{0}:
 &\ell^{2}(G\times Z)&\to&\ell^{2}(G/H \times \br{G}/\br{H}\times Z)
\\ &g &\mapsto& \Big[(x', y', z)\mapsto \esp{x \in p^{-1}(x')}
\big(\kappa(x-s(x'),\sigma(y'))g(x,\cdot)\big)(z)\Big] \end{matrix}, $$ $$\begin{matrix} T_{1}:
 &\ell^{2}(G/H \times \br{G}/\br{H}\times Z)&\to&\ell^{2}(\br{G}\times Z)
\\ &g &\mapsto& \Big[(y, z)\mapsto \esp{x' \in G/H}
\big(\kappa(s(x'),y)g(x',\pi(y), \cdot)\big)(z) \Big] \end{matrix}.$$ Then we have $T_{1}
\circ T_{0}=T_{G,\br{G},Z},$ since $\kappa|_{H\times \bar
H}=\mathrm{Id}_{\ell^2(Z)}.$

Recall that $T_{H,\br{G}/\br{H},Z}$ and $T_{G/H,\br{H},Z}$ are as follows
$$\begin{matrix} T_{H,\br{G}/\br{H},Z}: &\ell^{2}(H\times Z) &\to & \ell^{2}(\br{G}/\br{H}\times Z)
\\ &g& \mapsto & \Big[(y',z)\mapsto\esp{x \in H}
\big(\kappa_{0}(x,y')g(x,\cdot)\big)(z)\Big],
\\T_{G/H,\br{H}, Z}:&\ell^{2}(G/H\times Z) &\to & \ell^{2}(\br{H}\times Z) \\
&g& \mapsto & \Big[(y,z)\mapsto\esp{x' \in G/H}
\big(\kappa_{1}(x',y)g(x',\cdot)\big)(z)\Big].\end{matrix}$$

As explained in paragraph 1, $T_{1}$ and $T_{2}$ are block diagonals and each
block is identified with $T_{H,\br{G}/\br{H}, Z}$ and $T_{G/H,\br{H},Z}$.
Precisely, for a fixed $x' \in G/H,$ we define the following
isometries
$$\begin{matrix} \alpha_{x'}:&
\ell^{2}(p^{-1}(x')\times Z) &\to& \ell^{2}(H\times Z)
\\ &g &\mapsto&
\big[(x,z) \mapsto g(x+s(x'),z)\big] \end{matrix}$$
$$\begin{matrix} \beta:&
\ell^{2}(\bar G/\bar H\times Z) &\to& \ell^{2}(G/H\times \bar G/\bar H\times Z)
\\ &g &\mapsto&
\big[(x'',y'',z) \mapsto g(y'',z)\big]. \end{matrix}$$
By imbedding $\ell^2(p^{-1}(x')\times Z)$ into $\ell^2(G\times Z)$ (as
subspace of functions with support in $p^{-1}(x')\times Z$), we have
$$T_0|_{\ell^2(p^{-1}(x')\times Z)}=\beta\circ T_{H,\bar G/\bar H,Z}\circ\alpha_{x'}\in\L(\ell^2(p^{-1}(x')\times Z),\ell^{2}(G/H\times\bar G/\bar H\times Z)).$$ Moreover since
$T_0(\ell^2(p^{-1}(x')\times Z))\subset\ell^2(\{x'\}\times \bar G/\bar H\times
Z)$ and $\alpha_{x'}\otimes 1_E$ and $\beta\otimes 1_E$ are isometries, we have
$$\|T_0\otimes 1_E\|_{\L(\ell^{2}(G\times Z,E),\ell^{2}(G/H \times \br{G}/\br{H}\times Z,E))}$$ $$\leq \max_{x'\in G/H}\|T_0|_{\ell^2(p^{-1}(x')\times Z)}\otimes
1_E\|_{\L(\ell^2(p^{-1}(x')\times Z,E),\ell^{2}(\{x'\}\times \bar G/\bar H\times
Z,E))}$$ $$\leq \|T_{H,\bar G/\bar H,Z}\otimes 1_E\|_{\L(\ell^2(H\times
Z,E),\ell^2(\bar G/\bar H\times Z,E)}.$$ Similarly for $T_1,$ we fix $y'\in\bar
G/\bar H$ and define $$\begin{matrix}
\gamma_{y'}:& \ell^{2}(\br{H}\times Z) &\to& \ell^{2}(\pi^{-1}(y')\times Z)
\\ &g &\mapsto& \big[(y,z) \mapsto g(y-\sigma(y'),z)\big] \end{matrix}$$
$$\begin{matrix}
 \delta_{y'}:& \ell^{2}(G/H\times\{y'\}\times Z) &\to& \ell^{2}(G/H\times
 Z) \\ &g &\mapsto&
 \Big[(x'',z)\mapsto\Big(\kappa(s(x''),\sigma(y'))g(x'',y',\cdot)\Big)(z)
 \Big]. \end{matrix}$$ We imbed $\ell^2(G/H\times\{y'\}\times Z)$ into
 $\ell^2(G/H\times\bar G/\bar H\times Z),$ and then have
 $$T_1|_{\ell^2(G/H\times\{y'\}\times Z)}=\gamma_{y'}\circ T_{G/H,\bar
 H,Z}\circ\delta_{y'}\in\L(\ell^2(G/H\times\{y'\}\times
 Z),\ell^2(\pi^{-1}(y')\times Z)).$$ Since $\|\gamma_{y'}\otimes 1_E\|= 1$
 and $\|\delta_{y'}\otimes 1_E\|\leq 1$  (which
follows from the assumption $\|\kappa(x,y)\otimes 1_E\|\leq 1$), we conclude
that $$\|T_1\otimes 1_E\|_{\L(\ell^2(G/H\times\bar G/\bar H\times
Z,E),\ell^2(\bar G\times Z,E))}\leq \|T_{G/H,\bar H,Z}\otimes
1_E\|_{\L\big(\ell^2(G/H\times Z,E),\ell^2(\bar H\times Z,E)\big)},$$ and then
statement (i) is proved when $n=2.$

Now begin the induction. Suppose the proposition holds for $n=m \geq 2$. By hypothesis we are given the following chains of abelian groups,
$$0=G_{0} \subset G_{1} \subset G_{2} \subset ... \subset G_{m+1},$$
$$0=\br{G}^{m+1} \subset \br{G}^{m} \subset \br{G}^{m-1} \subset ... \subset \br{G}^{0}.$$
We apply the induction hypothesis to the following chains of abelian groups
$$0=G_{1}/G_{1} \subset G_{2}/G_{1} \subset ... \subset G_{m+1}/G_{1},$$
$$0=\br{G}^{m+1} \subset \br{G}^{m} \subset ... \subset \br{G}^{1},$$ and to
the map $\kappa':G_{m+1}/G_1\times \bar G^1\to\mathcal U(\ell^2(Z))$
defined by $\kappa'([x],y)=\kappa(x,y)$ for any $x\in G_{m+1},y\in \bar G^1,$
where $[\cdot]:G_{m+1}\to G_{m+1}/G_1$ is the projection. Then we get
$$\|T_{G_{m+1}/G_{1},\br{G}^{1},Z} \otimes 1_{E}\| \leq \prod_{i=1}^{m} \|T_{G_{i+1}/G_{i},\br{G}^{i}/\br{G}^{i+1},Z} \otimes 1_{E}\|.$$ We apply the proposition in the case when $n=2$, which is already proved, to
$$0=G_{0} \subset G_{1} \subset G_{m+1},$$
$$0=\br{G}^{m+1} \subset \br{G}^{1} \subset \br{G}^{0},$$
and get
$$\|T_{G_{m+1},\br{G}^{0},Z} \otimes 1_{E}\| \leq \|T_{G_{m+1}/G_{1},\br{G}^{1},Z} \otimes 1_{E}\| \cdot \|T_{G_{1},\br{G}^{0}/\br{G}^{1},Z} \otimes 1_{E}\|.$$
Combining this inequality with the previous one, we complete the proof of statement (i).

To prove (ii) we first show that for any $i\in I,$ $\br{G^i}/\br{G^{i+1}}$ is in
bijection with $ \widehat{G_{i+1}/G_i}.$ In fact, we define
 the following map $$\begin{matrix} \lambda^{*}_i: &\br{G^i}/\br{G^{i+1}} &\to&
 \widehat{G_{i+1}/G_i} \\ &y &\mapsto & \Big[x\mapsto \lambda_i(x,y)\Big]    
 \end{matrix}.$$ The non degeneracy of $\lambda_i(\cdot,y)$ implies injectivity,
 hence $$|G_{i+1}/G_i|=|\widehat{G_{i+1}/G_i}|\geq|\br{G^i}/\br{G^{i+1}}|.$$ By
 symmetry of $G$ and $\br{G}$ we have
 $|G_{i+1}/G_i|\leq|\br{G^i}/\br{G^{i+1}}|$, and hence $\lambda^{*}_i$ is a
 bijection.
 
 By identifying $\br{G^i}/\br{G^{i+1}}$ with $ \widehat{G_{i+1}/G_i},$
 we have $$T_{G_{i+1}/G_i,\bar
 G^i/\bar
 G^{i+1},Z}=T_{G_{i+1}/G_{i}}\otimes
 \mathrm{Id}_{\ell^2(Z)}\in\L(\ell^2(G_{i+1}/G_i\times Z),\ell^2(\widehat{G_{i+1}/G_i}\times Z)),$$ and thus $$\|T_{G_{i+1}/G_i,\bar
 G^i/\bar
 G^{i+1},Z}\otimes 1_E\|_{\L(\ell^2(G_{i+1}/G_i\times Z,E),\ell^2(\bar G^i/\bar
 G^{i+1}\times Z,E))}$$$$=\|T_{G_{i+1}/G_i}\otimes
 1_E\|_{\L(\ell^2(G_{i+1}/G_i,E),\ell^2(\widehat{G_{i+1}/G_{i}},E))},$$ which
 together with statement (i) implies statement (ii) immediately.
{\cqfd}

Next we want to prove the following corollary $\ref{VFFT2}$, which is the same as Lemme 3.2 in $\cite{JTA}$.

 Let $F$ be a non archimedian local field, $\O$ its ring of integers and $\pi$ a uniformizer of $\O$. Let $\E$ be a class of Banach spaces, $h \in \N$ and $\alpha>0$ such that for all $E \in \E$ we have $\|T_{\OO{h}} \otimes 1_{E}\| \leq e^{-\alpha}$.

 Let $\chi$ a non degenerate character of $\OO{h}$, which means that the restriction of $\chi$ over $\pi^{h-1}\OO{h}$ is non trivial. We recall that as a consequence the group morphism
\begin{align*}
\OO{h} &\rightarrow \widehat{\OO{h}}
\\ y &\mapsto (x \mapsto \chi(xy))
\end{align*}
is an isomorphism.

For any finite set $I$, define $\ell^{2}(I \times \OO{n})_{\chi}$ the subspace
of $\ell^{2}(I \times \OO{n})$ consisting of functions satisfying
$f(x,t+\pi^{n-h}s)=f(x,t)\chi(s)$, for $x \in I, t \in \OO{n}$, and $s \in
\OO{h}$. Let $Z$ be any set of representatives of $\OO{n-h}$ in $\OO{n}$.
Identify $\ell^2(I\times \OO{n})_\chi$ to $\ell^2(I\times Z)$ by restriction,
and for any operator $A\in\mathcal{L}(\ell^2(I\times\OO{n})_\chi)$ and any
Banach space $E$, equip $A\otimes 1_E$ with the norm defined in paragraph one.
This norm is independent of the choice of $Z$. In fact, for another set of
representatives $Z',$ there exists a map $\alpha$ from $Z$ to $\OO{h}$ such
that $Z'=\{z+\pi^{n-h}h(z)|z\in Z\}.$ And for any $i,j\in I$ and $z,w\in Z,$ the
matrix coefficient of $A$ for $\delta_{i,z+\pi^{n-h}\alpha(z)}$ and $\delta_{j,w+\pi^{n-h}\alpha(w)}$ is
$\chi(\frac{\alpha(z)}{\alpha(w)})$ times that for $\delta_{i,z}$ and
$\delta_{j,w}.$

\begin{cor}\label{VFFT2}
For any $n \in \mathbb{N}$ such that $n \geq h$, the operator
\begin{align*}
T: \ell^{2}(\O / \pi^{n}\O \times \O / \pi^{n}\O)_\chi &\rightarrow \ell^{2}(\O / \pi^{n}\O \times \O / \pi^{n}\O)_\chi
\\ (\xi_{x,t})_{x,t \in \O /\pi^{n}\O} & \mapsto (\esp{x \in \O/\pi^{n}\O}\xi_{x,t+xy})_{y,t \in \O/\pi^{n}\O}
\end{align*}
satifies $\| T \otimes 1_{E} \| \leq e^{-(\frac{n}{h}-1)\alpha}$ for any $E \in \E$.
\end{cor}

In fact, let $Z$ be a set of representatives of $\OO{n-h}$ in $\OO{n}$. Define
bilinear map $\kappa$ from $\OO{n} \times \OO{n}$ to the set of unitary
operators on $\ell^2(\O / \pi^{n}\O)_\chi$ (and hence $\ell^2(Z)$), by
$\kappa(x,y)(f)(t)=f(t+xy),$ for any $x,y,t \in \OO{n}, f \in \ell^2(\O /
\pi^{n}\O)_\chi.$ Let $n=ah+r,0 \leq r <h$. Apply proposition \ref{VFFT} to
$Z,\kappa,T,$ the following chains of finite abelian groups of length $a+1$ 
$$0 \subset \pi^{n-h}\OO{n} \subset \pi^{n-2h}\OO{n} \subset ... \subset
\pi^{n-ah}\OO{n} \subset \OO{n},$$ $$0 \subset \pi^{ah}\OO{n} \subset
\pi^{(a-1)h}\OO{n} \subset ... \subset \pi^{h}\OO{n} \subset \OO{n},$$ and
$I=\{0,1,...,a-1\}\subset\{0,1,...,a\}$. Then by corollary \ref{OF-alpha} we
complete the proof. 


{\cqfd}

\section{Strong Banach property (T) for $SL_{3}$ over a non archimedian local field}
We follow paragraph 4 in \cite{JTA} and \cite{duke}.

Let $F$ be a non archimedian local field, $\z$ its ring of integers, $\p$ a uniformizer, $\F$ the residue field.

The following theorem implies immediately that $SL_{3}(F)$ has strong Banach property $(T)$ in the sense of definition \ref{renf-ban}.

We note $G=SL_{3}(F)$ and $K=SL_{3}(\z)$.

We endow $G$ with the length $\ell$ defined by
 $$\ell(k
(\p^{\frac{i+2j}{3}}\begin{pmatrix} \p^{-(i+j)} & 0 & 0\\
                              0      &  \p^{-j}&0\\
                               0  &        &1\end{pmatrix})
k')=i+j$$ for $k,k'\in K$ and $i,j\in \N$, $i-j\in 3\Z$.

We recall that if $\ell'$ is a length over $G$, and $\E$ is a class of Banach spaces stable under duality and complex conjugation, we note $\E_{G,\ell'}$ the class of continuous representations $(E,\pi)$ of $G$ in a Banach space $E$ of the class $\E$ such that $\|\pi(g)\|_{\L(E)}\leq e^{\ell' (g)}$. We recall that $\CC^{\E}_{\ell'}(G)$ is the Banach algebra, endowed with complex conjugation and the usual involution, by completion of $C_{c}(G)$ for the norm $\|f\|=\sup_{(E,\pi)\in \E_{G,\ell'}}\|\pi(f)\|_{\L(E)}$.

 Let $h\in \N^{*},\alpha>0$ and $\E$ be a class of Banach spaces stable under duality and complex conjugation, such that for any $E$ in $\E$, we have
\begin{equation}
\|T_{\O/\pi^{h}\O}\otimes 1_{E}\|\leq e^{-\alpha} \tag{**} \label{type}.
\end{equation}

\begin{thm}\label{thmnonarchbanach}
Let $\beta \in [0,\frac{\alpha}{3h})$. There exists $t,C'>0$ such that for any
$C \in \R_{+}$, there exists a real and self-adjoint idempotent element $\pro\in
\CC^{\E}_{C+\beta\ell}(G)$ such that \begin{itemize} 
	\item (i) for any representation $(E,\pi) \in \E_{G,C+\beta \ell }$, the image
	of $\pi(\pro)$ is the subspace of E consisting of all $G$-invariant vectors, 
	\item (ii) there exists a sequence $\pro_{n}\in C_{c}(G)$, such that
	$\int_{G}|\pro_{n}(g)|dg\leq 1$, $\pro_{n}$ has support in $\{g\in
	G,\ell(g)\leq n\}$, and $\|\pro-\pro_{n}\|_{\CC^{\E}_{C+\beta\ell}(G)}\leq C'
	e^{2C-tn}$.
\end{itemize}
\end{thm}

The rest of this paragraph is devoted to the proof of theorem \ref{thmnonarchbanach}. We fix $h,\alpha$ and $\E$ as above. We note $\Lambda=\{(i,j)\in \N^{2}, i-j\in 3\Z\}$.

\begin{prop}\label{prop-decr-ban}
If $\beta \in [0,\frac{\alpha}{3h})$. There exists $C'>0$ such that the following property is satisfied. Let $C\in \R^{+}$ be arbitrary. If $(E,\pi)$ is in the class $\E_{G,C+\beta \ell }$ and $\xi, \eta$ are two $K$-invariant vectors of norm $1$ in $E$ and $E^{*}$, and if we put $c(g)=\s{\eta,\pi(g)\xi}$ for $g\in G$, then the function $c:\Lambda\to \C$ defined by abuse of notation by $c(i,j)=c(\p^{\frac{i+2j}{3}}\begin{pmatrix} \p^{-(i+j)} & 0 & 0\\
                              0      &  \p^{-j}&0\\
                               0  &        &1\end{pmatrix})
$ tends to a limit $c_{\infty}$ at infinity and we have $|c(i,j)-c_\infty|\leq C'e^{2C-(\frac{\alpha }{3h}-\beta)(i+j)}$.
\end{prop}

\begin{prop}\label{prop-decr-ban-equiv}
Let $\beta \in [0,\frac{\alpha}{3h})$. Let $(V,\pi)$ be a non trivial unitary
 irreducible representation of $K$. There exists a constant $C'>0$ such that the
 following property is satisfied. Let $C\in \R^{+}$ be arbitrary. If $(E,\pi)$
 is in the class $\E_{G,C+\beta \ell }$ and $\xi$ is a $K$-invariant vector of
 norm $1$ in $E$, and $\eta$ a $K$-invariant vector of norm $1$ in $V \otimes
 E^{*}$, and if we put $c(g)=\s{\eta,\pi(g)\xi}\in V$ for $g\in G$, then the
 function $c:\Lambda\to V$ defined by abuse of notation by
 $c(i,j)=c(\p^{\frac{i+2j}{3}}\begin{pmatrix} \p^{-(i+j)} & 0 & 0\\
                              0      &  \p^{-j}&0\\
                               0  &        &1\end{pmatrix})
$ tends to $0$ at infinity and we have $\|c(i,j)\|_{V}\leq C'e^{2C-(\frac{\alpha }{3h}-\beta)(i+j)}$.
\end{prop}

Propositions \ref{prop-decr-ban} and \ref{prop-decr-ban-equiv} imply theorem \ref{thmnonarchbanach}. It is formally identical to the argument in \cite{duke} showing that propositions 3.3 and 3.4 imply theorem 3.2 (this argument is placed between the statement of propositions 3.4 and 3.6 in \cite{duke}).

We now show that propositions \ref{prop-decr-ban} and \ref{prop-decr-ban-equiv} imply theorem \ref{thmnonarchbanach}. Let $\alpha >0$ and $\beta \in [0,\frac{\alpha}{3h})$. Proposition \ref{prop-decr-ban} allows the construction of an idempotent $\pro$ which thanks to proposition \ref{prop-decr-ban-equiv} is the one in theorem \ref{thmnonarchbanach}.

 Let $C \in R_{+}$. For $g \in G$ we put $P_{g}=e_{K}e_{g}e_{K}$ (which is of
 integral $1$ since $$P_g=\mathrm{vol}(K\cap
 gKg^{-1})\chi_{KgK}=\mathrm{vol}(\mathrm{Stab}_K(gK))\chi_{KgK},$$ and
 $\mathrm{vol}(\mathrm{Stab}_K(gK)\mathrm{vol}(KgK))=\mathrm{vol}(K)=1$).
 Following from proposition \ref{prop-decr-ban}, $P_{g}$ is a Cauchy sequence in $\CC^{\E}_{C+\beta \ell}(G)$ when $g$ tends to infinity, and this sequence
 tends to a certain element $\pro$. Since $P^{*}_{g}=P_{g^{-1}}$, $\pro$ is
 self-adjoint. Furthermore for any representation $(E,\pi)$ in the class
 $\E_{G,C+\beta\ell}$, the image of $\pi(\pro)$ consists of elements $x$ such
 that for all $g \in G$, $\pi(e_{K})\pi(g)x=x$. In fact $e_{K}e_gP_{g'}$ is
 equal to $\int_{K}P_{gkg'} dk$, thus it follows that $e_{K}e_{g}\pro=\pro$. It
 is also evident that the restriction of $\pi(\pro)$ on the space of vectors $x$
 such that for all $g \in G$, $\pi(e_{K})\pi(g)x=x$ is $\mathrm{Id}$, thus
 $\pro$ is an idempotent. If $\pi$ is isometric and $E$ is uniformly convex, the
 following lemma shows that the vectors $x$ such that for all $g \in G$,
 $\pi(e_{K})\pi(g)x=x$ are exactly $G$-invariant vectors.                            

\begin{lem}\label{cas-isom}
Let $\pi$ be an isometric representation of $G$ over a uniformly convex space $E$. Let $x\in E$ be of norm $1$, such that $\pi(e_{K})\pi(g)x=x$ for all $g \in G$. Then $x$ is $G$-invariant.
\end{lem}

In fact, denoting $f(k)=\pi(k)\pi(g)x,\forall k \in K$, we have $\|f(k)\|_{E}=1$ since $\pi$ is isometric, and $\|\int_{K}{f(k)d k}\|_{E}=\|x\|_{E}=1$ since $f:K \to E$ is continuous. Since $E$ is uniformly convex, $f(k)$ is constant. Thus we have $x=\pi(g)x$.


\cqfd

In general we use proposition \ref{prop-decr-ban-equiv}. Let $(V,\rho)$ be a non
trivial irreducible unitary representation of $K$. It is of finite dimension
since $K$ is compact. Let $e_{K}^{V}=\int_{K}\chi_{V}(k)^{*}e_{k}dk$, where
$\chi_{V}$ is the character of $(V,\rho)$ (i.e.,
$\chi_{V}(k)=\sum_{i=1}^{\text{dim}V}\s{e_{i},\rho(k)e_{i}}$ for a basis
$\{e_{1},...,e_{\text{dim}V}\}$ of $V$). Proposition \ref{prop-decr-ban-equiv}
implies that $e_{K}^{V}e_{g}e_{K}$ tends to $0$ in $\CC^{\E}_{G,C+\beta\ell}$ when $g$
tends to infinity in $G$. It follows that $e_{K}^{V}e_{g}\pro=0$ in
$\CC^{\E}_{G,C+\beta\ell}$, and that the image of $\pi(\pro)$ consists of
exactly $G$-invariant vectors of $H$. Finally we put $\pro_{n}=P_{g_{E(n/2)}}$,
with $g_{n}=\p^{n}\begin{pmatrix} \p^{-2n} & 0 & 0\\
                              0      &  \p^{-n}&0\\
                               0  &        &1\end{pmatrix}$. We verify easily that $t$ and $C'$ are independent of $C$. This completes the proof of theorem \ref{thmnonarchbanach}.

                                \cqfd



Let $G'=\{g\in GL_{3}(F),\det g\in \p^{\Z}\}/(\p^{\Z}\mathrm{Id})$. Then $G$ is a subgroup of index $3$ in $G'$ and $K$ is a maximal compact subgroup of $G'$. Any element of $G'$ can be written in the form $k\begin{pmatrix} \p^{-(i+j)} & 0 & 0\\
                              0      &  \p^{-j}&0\\
                               0  &        &1\end{pmatrix}
                               k'$, with unique $i,j\in \N$ and $k,k'\in K$. We endow $G'$ of the length defined by  $\ell(k
\begin{pmatrix} \p^{-(i+j)} & 0 & 0\\
                              0      &  \p^{-j}&0\\
                               0  &        &1\end{pmatrix}
k')=i+j$, for $k,k'\in K$ and $i,j\in \N$.

By inducing the representation of $G$ to $G'$, we see that the propositions \ref{prop-decr-ban} and \ref{prop-decr-ban-equiv} follow from the following two propositions.


Precisely, let $(E,\pi)$ be a continuous representation of $G$, and
$\text{Ind}_{G}^{G'}E$ be the space of continuous maps $f$ from $G'$ to $E$ such
that $f(xg)=\pi(g^{-1})f(x)$ for all $g\in G, x\in G'$, endowed with the norm
given by $$\|f\|:=\sum_{[x]\in G'/G}\|f(s([x]))\|_{E},$$where $s:G'/G \simeq
\Z/3\Z \to G'$ is a section such that $s([e])=e$. By simple calculation one
verifies that the following map is isometric imbedding from $\pi$ to
$\text{Ind}_{G}^{G'}\pi|_{G}$ $$\begin{matrix} i:&E &\to&
\text{Ind}_{G}^{G'}E& \\ &\xi &\mapsto& \Big(x\mapsto \Big\{\begin{matrix}0&
\text{ if } x\not\in G \\ \pi(x^{-1})\xi& \text{ if } x\in G \end{matrix}\Big)&,
\end{matrix}$$
and that $\s{\eta,\xi}=\s{e^{*}\eta,i\xi}$ for any $\eta \in E^{*}$ and $\xi \in E$, where $e^{*}$ is the adjoint of evaluation at identity $e \in G':\text{Ind}_{G}^{G'}E\to E$. One sees that propositions \ref{prop-decr-ban} and \ref{prop-decr-ban-equiv} follow immediately from the following two propositions.

\begin{prop}\label{prop-decr-ban'}

 Let $\alpha>0$ and $\beta\in [0,\frac{\alpha }{3h})$. There exists $C'>0$ such that the following properties are satisfied. Let $C \in \R^{+}$ be arbitrary. If $(E,\pi)$ is in the class $\E_{G',C+\beta\ell}$ and $\xi$ and $\eta$ are two $K$-invariant vectors of norm $1$ in $E$ and $E^{*}$, and if we put $c(g)=\s{\eta,\pi(g)\xi}$ for $g\in G'$, then the restriction to $\Lambda$ of the function $c:\N^{2}\to \C$ defined by abuse of notation by $c(i,j)=c(\begin{pmatrix} \p^{-(i+j)} & 0 & 0\\
                              0      &  \p^{-j}&0\\
                               0  &        &1\end{pmatrix})
$ tends to a limit $c_{\infty}$ at infinity and we have $|c(i,j)-c_\infty|\leq C'e^{2C-(\frac{\alpha }{3h}-\beta)(i+j)}$ for $(i,j)\in \Lambda$.
\end{prop}

\noindent {\bf Remark.} By essentially the same method (as explained below) one can show that proposition \ref{prop-decr-ban'} still holds for $\{(i,j)\in\N^{2},i-j\in3\Z+1\}$ and $\{(i,j)\in\N^{2},i-j\in3\Z+2\}$, but the limits might be different from $c_{\infty}$. Indeed this happens when $\pi$ is a non trivial character of $G'/G=\Z/3\Z$.

\begin{prop}\label{prop-decr-ban-equiv'}
 Let $\alpha>0$ and $\beta\in [0,\frac{\alpha }{3h})$, and $(V,\tau)$ a non trivial unitary irreducible representation of $K$. There exists $C'>0$ such that the following properties are satisfied. Let $C \in \R^{+}$ be arbitrary. If $(E,\pi)$ is in the class $\E_{G',C+\beta\ell}$ and $\xi$ is a $K$-invariant vector of norm $1$ in $E$, and $\eta$ a $K$-invariant vector of norm $1$ in $V \otimes E^{*}$, and if we put $c(g)=\s{\eta,\pi(g)\xi} \in V$ for $g\in G'$, then the restriction to $\Lambda$ of the function $c:\N^{2}\to V$ defined by abuse of notation by $c(i,j)=c(\begin{pmatrix} \p^{-(i+j)} & 0 & 0\\
                              0      &  \p^{-j}&0\\
                               0  &        &1\end{pmatrix})
$ tends to $0$ at infinity and we have $\|c(i,j)\|_{V}\leq C'e^{2C-(\frac{\alpha }{3h}-\beta)(i+j)}$ for $(i,j)\in \Lambda$.

\end{prop}

It remains to prove propositions \ref{prop-decr-ban'} and
\ref{prop-decr-ban-equiv'}. We begin by the following preliminary lemma. 


\begin{lem}\label{compen}
 Let $\chi:\F\to \C^{*}$ be a non trivial character. Let $h\in \N^{*}, \alpha\in \R_{+}^{*}, n\in \N^{*}$. Let $E$ be a Banach space such that $\|T_{\O/\pi^{h}\O}\otimes 1_{E}\|\leq e^{-\alpha}$, and $(\xi_{x,y})_{x,y\in \z/\p^{n}\z}$ vectors of $E$. Then
$$\esp{a,b\in \z/\p^{n}\z}
\Big\|\esp{x\in \z/\p^{n}\z,\eps \in \F} \chi(\eps)\xi_{x,ax+b+\p^{n-1}\eps}\Big\|^{2} $$ $$\leq
(\sharp \O/\p^{h-1}\O)^{2}
e^{-2(\frac{n}{h}-1)\alpha}\esp{x,y\in \z/\p^{n}\z}
\|\xi_{x,y}\|^{2}.$$
\end{lem}

In fact, the statement is obvious for $n<h$. We suppose then $n \geq h$. We have
$$\esp{x\in \z/\p^{n}\z,\eps \in \F} \chi(\eps)\xi_{x,ax+b+\p^{n-1}\eps}=\sum_{\eta} \Big(\esp{x\in \z/\p^{n}\z,z \in \O/\pi^{h}\O} \eta(z)\xi_{x,ax+b+\p^{n-h}z}\Big)$$
where the sum is taken over characters $\eta$ of $\OO{h}$ whose restriction to $\pi^{h-1}\O/\pi^{h}\O=\O/\pi\O=\F$ is equal to $\chi$. Such a character $\eta$ is non degenerate. It then suffices to show that for any non degenerate character $\eta$ of $\OO{h}$ and for any vectors $(\xi_{x,y})_{x,y\in \z/\p^{n}\z}$ in $E$ we have
$$\esp{a,b\in \z/\p^{n}\z}
\Big\|\esp{x\in \z/\p^{n}\z,z\in \O/\pi^{h}\O}\eta(z)\xi_{x,ax+b+\p^{n-h}z}\Big\|^{2} $$ $$\leq
e^{-2(\frac{n}{h}-1)\alpha}\esp{x,y\in \z/\p^{n}\z}
\|\xi_{x,y}\|^{2}.$$
This follows immediately from corollary \ref{VFFT2} by applying it to
$$(\esp{z\in\OO{h}}\eta(z)\xi_{x,y+\pi^{n-h}z})_{x,y\in\OO{n}}.$$ {\cqfd}

Afterwards we note $B$ the building of $PGL_{3}(F)$. We recall that the vertices of $B$ are identified with lattices of $F^{3}$, up to dilations, and that $PGL_{3}(F)$ acts on the left over $B$. Thus $G'$ acts on the left over $B$ and this action is transitive. On the other hand given $x_{0}$ a vertex corresponding to lattice $\O^{3}$ in $F^{3}$, its stabilizer in $G'$ is $K=SL_{3}(\mathcal{O})$, which allows the identification of $B$ and $G'/K$. If $M$ is a lattice associated to a vertex $x$ of $B$, $\det(M)=\p^{-a}\det(\z^{3})$ for a certain integer $a \in \Z$, whose image in $\Z/3\Z$ does not depend on the choice of representatives of $x$ , we say that it is the type of $x$. Given $x,y \in B$ there exist unique integers $i,j \in \N$ such that for a certain basis $(v_1,v_2,v_3)$ of $F^{3}$ generating $x$ as $\mathcal{O}$-module we have $y=\z \p^{-i-j}v_1+\z \p^{-j}v_2+\z v_3 \text{ modulo }F^{*}$. We write $\sigma(x,y)=(i,j) \in \N^{2}$. Then we have $\sigma(y,x)=(j,i)$ and the type of $y$ is type($x$)$+i-j$ modulo $3$. Now given $x,y\in B=G'/K$, we have $\sigma(x,y)=(i,j)$ if and only if
$x^{-1}y=K
\begin{pmatrix} \p^{-(i+j)}&0&0\\0&\p^{-j}&0\\0&0&1\end{pmatrix}
K$ in $K\backslash G'/K$.

We note $(e_{1},e_{2},e_{3})$ the base of $F^{3}$. Let $m,n \in \N$. For $x,y\in \z/\p^{n}\z$ and $a,b\in \z/\p^{m}\z$ we note $M^{n}_{x,y}$  and $M^{-m}_{a,b}$ the following lattices of $F^3$:
$$M^{n}_{x,y}=\z\p^{-n}(e_{1}+xe_{2}-ye_{3})+\z e_{2}+\z e_{3}\text{ and}$$
$$M^{-m}_{a,b}=\{u_{1}e_{1}+u_{2}e_{2}+u_{3}e_{3}\in \z^{3},
u_{1}b+u_{2}a+u_{3}\in \p^{m}\z\}$$ $$=
\z(e_{1}-be_{3})+\z (e_{2}-ae_{3})+\z \p^{m}e_{3}.$$
By abuse of notation we have supposed that $x,y,a,b$ were lifted to be elements of $\O$, and the result does not depend on such choices. We note still $M^{n}_{x,y}$  and $M^{-m}_{a,b}$ the dilation classes of these lattices, considered as elements of $B=G'/K$. In $G'/K$ we have
$M^{n}_{x,y}= \begin{pmatrix} \p^{-n}&0&0\\ \p^{-n}x&1&0\\ -\p^{-n}y&0&1
\end{pmatrix}    K$
and
$M^{-m}_{a,b}=    \begin{pmatrix} 1&0&0\\ 0&1&0 \\ \p^{-m}b&\p^{-m}a&\p^{-m}
\end{pmatrix}^{-1}  K$.

\noindent {\bf Remark} $M_{a,b}^{-m}$ and $M_{x,y}^{n}$ belong to the $K$-orbit of $\begin{pmatrix} 1&0&0\\ 0&1&0\\ 0&0&\p^{m}
\end{pmatrix}    K$ and $\begin{pmatrix} \p^{-n}&0&0\\ 0&1&0\\ 0&0&1
\end{pmatrix}    K$ in $G'/K$, which are identified with the projective planes $\mathbb{P}^{2}(\OO{m})$ and $\mathbb{P}^{2}(\OO{n})$. In fact, the points $M_{a,b}^{-m}$ for $a,b \in \OO{m}$ and $M^{n}_{x,y}$ for $x,y\in \z/\p^{n}\z$ describe the affine planes in these projective planes.

The following lemma shows that the relative position of $M^{-m}_{a,b}$ and
$M^{n}_{x,y}$ only depends on the scalar product between the vector
$e_{1}+xe_{2}-ye_{3}$ and the linear form
$u_{1}e_{1}+u_{2}e_{2}+u_{3}e_{3}\mapsto u_{1}b+u_{2}a+u_{3}$ (we choose these
letters since it takes the agreeable form $(ax+b)-y$).     

\begin{lem}\label{6}

 Let $x,y\in \z/\p^{n}\z$ and $a,b\in \z/\p^{m}\z$. Let $i$ be the largest integer of $\{0,...,\min(m,n)\}$ such that $y-(ax+b)$ belongs to $\p^{i}\z/\p^{\min(m,n)}\z$. Then $\sigma(M^{-m}_{a,b},M^{n}_{x,y})=(m+n-2i,i)$.
\end{lem}

In fact if $i<min(m,n)$ we have
$$M^{n}_{x,y}=\z \p^{-n}(e_{1}+xe_{2}-ye_{3})+\z (e_{2}-ae_{3})+\z \p^{-i}(e_{1}+xe_{2}-(ax+b)e_{3})\text{ and }$$ $$
M^{-m}_{a,b}=\z \p^{m-i}(e_{1}+xe_{2}-ye_{3})+\z (e_{2}-ae_{3})+\z (e_{1}+xe_{2}-(ax+b)e_{3})$$
if $i=n$ (and then $m\geq n$) we have
$$M^{n}_{x,y}=\z \p^{-n}(e_{1}+xe_{2}-(ax+b)e_{3})+\z (e_{2}-ae_{3})+\z e_{3}\text{ and } $$ $$
M^{-m}_{a,b}=\z (e_{1}+xe_{2}-(ax+b)e_{3})+\z (e_{2}-ae_{3})+\z \p^{m}e_{3}$$
and if $i=m$ (and then $n\geq m$) we have
$$M^{n}_{x,y}=\z \p^{-n}(e_{1}+xe_{2}-ye_{3})+\z (e_{2}-ae_{3})+\z e_{3}\text{ and } $$ $$
M^{-m}_{a,b}=\z (e_{1}+xe_{2}-ye_{3})+\z (e_{2}-ae_{3})+\z \p^{m}e_{3}.$$\cqfd

Now we prove proposition \ref{prop-decr-ban'}. Let $m,n \in \N$, with $m \geq n$. Let $x,y \in \OO{n}$ and $a,b \in \OO{m}$. We put $\xi_{x,y}=\pi(M^{n}_{x,y})\xi\in E$ and $\eta_{a,b}={}^{t}\pi((M^{-m}_{a,b})^{-1})\eta\in E^{*}$. Then $\|\xi_{x,y}\|\leq e^{C+n\beta}$, $\|\eta_{a,b}\|\leq e^{C+m\beta}$ and
$\s{\eta_{a,b},\xi_{x,y}}=c(m+n-2i,i)$, where $i$ is the largest integer of $\{0,...,n\}$ such that $y-(ax+b)$ belongs to $\unif^{i}\OO{n}$. It follows from Cauchy-Schwarz inequality and lemma \ref{compen} that
$$\Big|
\esp{a,b\in  \z/\p^{m}\z,x\in \z/\p^{n}\z,\eps \in \F}
\chi(\eps)\s{\eta_{a,b},\xi_{x,ax+b+\eps \p^{n-1}}}\Big|$$
$$\leq \sqrt{ \esp{a,b\in  \z/\p^{m}\z}
\|\eta_{a,b}\|^{2}}$$ $$
\sqrt{\esp{a,b \in \z/\p^{m}\z}
\Big\|\esp{x\in \z/\p^{n}\z,\eps \in \F}
\chi(\eps)\xi_{x,ax+b+\eps \p^{n-1}}\Big\|^{2}}$$
$$\leq q^{h-1}e^{-(\frac{n}{h}-1)\alpha}e^{2C+(m+n)\beta}.$$

But the left hand side is equal to
$\frac{1}{q}|c(m-n,n)-c(m-n+2,n-1)|$,
and then $|c(m-n,n)-c(m-n+2,n-1)|\leq q^{h}e^{-(\frac{n}{h}-1)\alpha}e^{2C+(m+n)\beta}$.

We have then, for $i,j \in \N$ with $j>0$,
$$|c(i,j)-c(i+2,j-1)|\leq q^{h}e^{-(\frac{j}{h}-1)\alpha}e^{2C+(i+2j)\beta}.
$$

By the automorphism

 $g\mapsto
\begin{pmatrix} 0&0&1\\ 0&1&0\\ 1&0&0\end{pmatrix}
{}^{t} g^{-1}\begin{pmatrix} 0&0&1\\ 0&1&0\\ 1&0&0\end{pmatrix}$ of $G'$, which stabilizes $K$, and maps $\begin{pmatrix} \p^{-(i+j)}&0&0\\0&\p^{-j}&0\\0&0&1\end{pmatrix}$ to $\begin{pmatrix} \p^{-(i+j)}&0&0\\0&\p^{-i}&0\\0&0&1\end{pmatrix}$, we obtain, for $i,j \in \N$ with $i>0$, $$|c(i,j)-c(i-1,j+2)|\leq q^{h}e^{-(\frac{i}{h}-1)\alpha}e^{2C+(2i+j)\beta}.$$.

With these two steps, we can obtain the limit $c_{\infty}$ by first moving from $(i,j)$ to the diagonal, and then along the diagonal to infinity as illustrated in the following Weyl chamber

\vskip 3mm

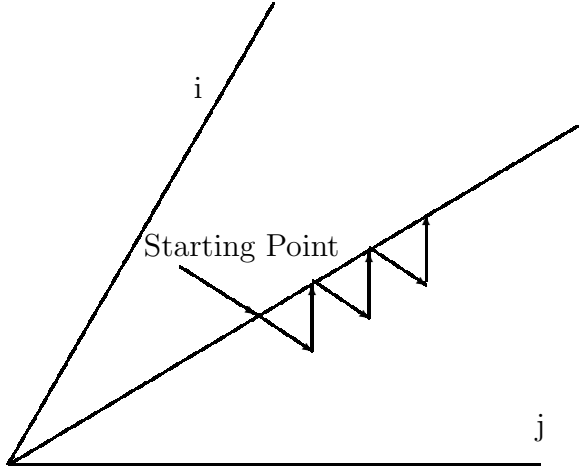
\ifx\JPicScale\undefined\def\JPicScale{1}\fi
\unitlength \JPicScale mm
\begin{picture}(120,60)(0,0)
\linethickness{0.3mm}
\put(10,0){\line(1,0){70}}
\linethickness{0.3mm}
\multiput(10,0)(0.12,0.21){292}{\line(0,1){0.21}}
\linethickness{0.3mm}
\multiput(10,0)(0.2,0.12){375}{\line(1,0){0.2}}
\put(120,30){\makebox(0,0)[cc]{}}

\put(35,50){\makebox(0,0)[cc]{i}}

\put(80,5){\makebox(0,0)[cc]{j}}

\linethickness{0.3mm}
\multiput(42.5,20)(0.18,-0.12){42}{\line(1,0){0.18}}
\put(50,15){\vector(3,-2){0.12}}
\linethickness{0.3mm}
\put(50,15){\line(0,1){8.75}}
\put(50,23.75){\vector(0,1){0.12}}
\put(40.75,28.75){\makebox(0,0)[cc]{Starting Point}}

\linethickness{0.3mm}
\linethickness{0.3mm}
\multiput(32.5,26.25)(0.19,-0.12){52}{\line(1,0){0.19}}
\put(42.5,20){\vector(3,-2){0.12}}
\linethickness{0.3mm}
\multiput(50,24.38)(0.18,-0.12){42}{\line(1,0){0.18}}
\put(57.5,19.38){\vector(3,-2){0.12}}
\linethickness{0.3mm}
\put(57.5,19.38){\line(0,1){8.75}}
\put(57.5,28.12){\vector(0,1){0.12}}
\linethickness{0.3mm}
\multiput(57.5,28.75)(0.18,-0.12){42}{\line(1,0){0.18}}
\put(65,23.75){\vector(3,-2){0.12}}
\linethickness{0.3mm}
\put(65,24.38){\line(0,1){8.75}}
\put(65,33.12){\vector(0,1){0.12}}
\end{picture}

\vskip 10mm

Precisely for $i,j \in \N$, with $i\geq j$ and $i-j \in 3\Z$, we have
$$|c(i,j)-c(\frac{2i+j}{3},\frac{2i+j}{3})|\leq
q^{h}e^{2C+\alpha}(e^{-i\frac{\alpha}{h}}+...+e^{-(\frac{2i+j}{3}+1)\frac{\alpha}{h}})e^{(2i+j)\beta}$$ $$
\leq q^{h}e^{2C+\alpha}(e^{\frac{\alpha}{h}}-1)^{-1}e^{-(2i+j)(\frac{\alpha}{3h}-\beta)},$$
and by the inequality of reversing $i$ and $j$, we have
$$|c(i,i)-c(i+1,i+1)|\leq
q^{h}e^{2C+\alpha}(1+e^{3\beta})e^{-i(\frac{\alpha}{h}-3\beta)}.$$

Proposition \ref{prop-decr-ban'} is then proved.\cqfd

Now we prove proposition \ref{prop-decr-ban-equiv'}(we recall that by lemma
\ref{cas-isom} this proposition is not necessary if we limit ourselves to
isometric actions over uniformly convex Banach spaces).     

For $x,y,a,b \in \O$, we must calculate the image of
$$\begin{pmatrix} 1&0&0\\ 0&1&0 \\ \p^{-m}b&\p^{-m}a&\p^{-m}
\end{pmatrix}
\begin{pmatrix} \p^{-n}&0&0\\ \p^{-n}x&1&0\\ -\p^{-n}y&0&1
\end{pmatrix}$$
in $K' \backslash G'/K$, where $K'$ is the open compact subgroup of $K$ included
in the kernel of $\tau$, then lemma \ref{6} calculate the image in $K \backslash
G'/K$. For simplicity we will only conduct the calculation when $m \geq n$ and
$y=ax+b$ or $y=ax+b+\p^{n-1}$. The lemma below, which is a variant of lemma
\ref{compen}, allows us to limit ourselves to these cases. The role of the
supplementary variable $k$ will be explicit later on.        

We know that any non trivial character of $\F$ is of the form $\chi_{d}:x \mapsto \chi(dx)$, for $d\in F^*$. We write  $\delta_{0}-\delta_{1}=\esp{d\in \F^{*}}t_{d}\chi_{d}$, for some certain $t_{d}\in \C$, and put $C_{2}=\esp{d\in \F^{*}}|t_{d}|^{2}$.

   \begin{lem}\label{compensanschi}
   Let $n\in \N^{*}$ and $k\in \{0,...,n\}$. Let $E$ be a Banach space in the class $\E$, $(\xi_{x,y})_{x,y\in \z/\p^{n}\z}$ vectors in $E$. Then
$$\esp{a,b\in \z/\p^{n}\z}
\Big\|\esp{x\in \p^{k}\z/\p^{n}\z} \xi_{x,ax+b}
-\esp{x\in \p^{k}\z/\p^{n}\z} \xi_{x,ax+b+\p^{n-1}}
\Big\|^{2}$$ $$\leq C_{2}q^{2h}e^{-2(\frac{n-k}{h}-1)\alpha}\esp{x\in  \p^{k}\z/\p^{n}\z,y\in \z/\p^{n}\z}
\|\xi_{x,y}\|^{2}.$$
\end{lem}

In fact, when $k=0$, we write
$$   \esp{x\in \z/\p^{n}\z} \xi_{x,ax+b}
-\esp{x\in \z/\p^{n}\z} \xi_{x,ax+b+\p^{n-1}}$$ $$=q\esp{d\in \F^{*}}t_{d}(
\esp{x\in \z/\p^{n}\z,\eps \in \F} \chi_{d}(\eps)\xi_{x,ax+b+\p^{n-1}\eps}),$$
then we apply Cauchy-Schwarz inequality and lemma \ref{compen}.

In general, note that if $a\in \OO{n},x\in \pi^{k}\OO{n}$, the product $ax\in \OO{n}$ only depends on $a$ modulo $\pi^{n-k}\O$. Applying the case when $k=0$ to $$(\xi_{\pi^{k}x_{1},\pi^{k}y_{1}+s([b])})_{x_{1},y_{1}\in \OO{n-k}},$$ for fixed $[b] \in \OO{k}$, where $s:\OO{k}\to \OO{n}$ is a section of $\OO{n}\to \OO{k}$, we then complete the proof.

\cqfd

\noindent{\bf Remark.}
1) By the same method we can get the same estimate for any $f\in \ell^{2}(\F)$ satisfying $\esp{\epsilon \in \F}f(\epsilon)=0$.

\noindent 2) Condition (\ref{type})
 (which implies $\E$ is of type $>1$) is necessary because otherwise we take $E=\ell^{1}((\OO{n})^{2})$ and $\xi_{x,y}=\delta_{x,y}\in E$, then the inequality of lemma \ref{compensanschi} becomes
 $$2\leq C_{2}q^{2h}e^{-2(\frac{n-k}{h}-1)\alpha},$$
 and it implies $\alpha=0$.

We fix now $k \in \N$ such that $\tau$ factorizes through $SL_{3}(\z/\p^{k}\z)$
and we note $K'$ the kernel of $K\to SL_{3}(\z/\p^{k}\z)$. We will apply lemma
\ref{compensanschi} for this value of $k$.     

Let $m,n\in \N$, with $m\geq n$. For $a,b \in \z$, the class in $G'/K'$ of $ \begin{pmatrix} 1&0&0\\ 0&1&0 \\ \p^{-m}b&\p^{-m}a&\p^{-m}
\end{pmatrix}^{-1}$ depends only on $a,b$ modulo $\p^{m+k}\z$. For $x,y\in \z/\p^{n}\z$ and $a,b\in \z/\p^{m+k}\z$ we define
 $$\xi_{x,y}=\pi
\begin{pmatrix} \p^{-n}&0&0\\ \p^{-n}x&1&0\\ -\p^{-n}y&0&1
\end{pmatrix} \xi \in E\text{ et }$$
$$\eta_{a,b}=({}^{t}\pi
\begin{pmatrix} 1&0&0\\ 0&1&0 \\ \p^{-m}b&\p^{-m}a&\p^{-m}
\end{pmatrix}\otimes 1)\eta\in E^{*}\otimes V.$$
Then
$\|\xi_{x,y}\|\leq e^{C+n\beta}$, $\|\eta_{a,b}\|\leq e^{C+m\beta}$.
We have $\s{\eta_{a,b},\xi_{x,y}}=c(A^{x,y}_{a,b})\in V$, where
$$A^{x,y}_{a,b}=
\begin{pmatrix} 1&0&0\\ 0&1&0 \\ \p^{-m}b&\p^{-m}a&\p^{-m}
\end{pmatrix}\begin{pmatrix} \p^{-n}&0&0\\ \p^{-n}x&1&0\\ -\p^{-n}y&0&1
\end{pmatrix}$$ $$=\begin{pmatrix} \p^{-n}&0&0\\ \p^{-n}x&1&0\\ \p^{-m-n}(ax+b-y)&\p^{-m}a&\p^{-m}
\end{pmatrix}.$$

We note $c(i,j)=c\begin{pmatrix} \p^{-(i+j)}&0&0\\ 0&\p^{-j}&0 \\ 0&0&1
\end{pmatrix}$. We recall that $c(kgk')=\tau(k)c(g)\in V$ for  $k,k'\in K,g\in G$. If $y=ax+b$ we have
$$A^{x,y}_{a,b}=\begin{pmatrix} -\p^{m-n+1}&1
&0\\ -\p^{m-n+1}x&x&1\\1&0&0
\end{pmatrix}
\begin{pmatrix} \p^{-m}&0&0\\ 0&\p^{-n}&0 \\ 0&0&1
\end{pmatrix}
\begin{pmatrix} 0&a&1\\ 1&\p a&\p \\ 0&1&0\end{pmatrix}$$
and then
$$c(A^{x,y}_{a,b})=\tau
\begin{pmatrix} -\p^{m-n+1}&1
&0\\ -\p^{m-n+1}x&x&1\\1&0&0
\end{pmatrix}
c\begin{pmatrix} \p^{-m}&0&0\\ 0&\p^{-n}&0 \\ 0&0&1
\end{pmatrix}$$ $$=\tau
\begin{pmatrix} -\p^{m-n+1}&1
&0\\ -\p^{m-n+1}x&x&1\\1&0&0
\end{pmatrix}c(m-n,n).$$

If $y=ax+b+\p^{n-1}$, we have
$$A^{x,y}_{a,b}=\begin{pmatrix} -\p^{m-n+1}&1&0\\ -\p^{m-n+1} x&x&1\\1&0&0
\end{pmatrix}\begin{pmatrix} \p^{-(m+1)}&0&0\\ 0&\p^{-(n-1)}&0 \\ 0&0&1
\end{pmatrix}\begin{pmatrix}
-1&\p a&\p \\0&a&1\\0&1&0
\end{pmatrix}$$
and then
$$c(A^{x,y}_{a,b})=\tau\begin{pmatrix} -\p^{m-n+1}&1&0\\-\p^{m-n+1} x&x&1\\1&0&0
\end{pmatrix}c\begin{pmatrix} \p^{-(m+1)}&0&0\\ 0&\p^{-(n-1)}&0 \\ 0&0&1
\end{pmatrix}$$ $$ =\tau\begin{pmatrix} -\p^{m-n+1}&1&0\\-\p^{m-n+1} x&x&1\\1&0&0
\end{pmatrix}c(m-n+2,n-1).$$

We recall that $\tau$ factorizes through $SL_{3}(\z/\p^{k}\z)$. We will apply
lemma \ref{compensanschi} to $k,$ then for any $x
\in\p^{k}\z/\p^{n}\z$,    
we always have $$\tau\begin{pmatrix} -\p^{m-n+1}&1&0\\-\p^{m-n+1} x&x&1\\1&0&0
\end{pmatrix}=\tau\begin{pmatrix} -\p^{m-n+1}&1&0\\0&0&1\\1&0&0
\end{pmatrix}$$ (this is the reason why we introduced the supplementary variable
$k$ in lemma \ref{compensanschi}).

By applying Cauchy-Schwarz inequality and lemma \ref{compensanschi}, we obtain
     $$\Big\|
\esp{a,b\in \z/\p^{m+k}\z,x\in \p^{k}\z/\p^{n}\z}
(\s{\eta_{a,b},\xi_{x,ax+b}}-
\s{\eta_{a,b},\xi_{x,ax+b+\p^{n-1}}})\Big\|_{V}$$
$$\leq \esp{a,b\in \z/\p^{m+k}\z}
\|\eta_{a,b}\| \Big\|
\esp{x\in \p^{k}\z/\p^{n}\z}
\xi_{x,ax+b}-\xi_{x,ax+b+\p^{n-1}}\Big\|$$
 $$\leq \sqrt{C_{2}} q^{h}e^{-(\frac{n-k}{h}-1)\alpha}e^{2C+(m+n)\beta}.$$

According to the preceding calculations, this inequality is rewritten as
 $$\Big\|\tau\begin{pmatrix} -\p^{m-n+1}&1&0\\0&0&1\\1&0&0
\end{pmatrix}c(m-n,n)-
\tau\begin{pmatrix} -\p^{m-n+1}&1&0\\0&0&1\\1&0&0
\end{pmatrix}c(m-n+2,n-1)\Big\|_{V}$$ $$\leq
      \sqrt{C_{2}} q^{h}e^{-(\frac{n-k}{h}-1)\alpha}e^{2C+(m+n)\beta},
      \text{ and then }$$
      $$\|c(m-n,n)-c(m-n+2,n-1)\|_{V}\leq
  \sqrt{C_{2}}
       q^{h}e^{-(\frac{n-k}{h}-1)\alpha}e^{2C+(m+n)\beta}.$$

We have then
$$\|c(0,3i)-c(2i,2i)\|_{V}\leq  \sqrt{C_{2}}
       q^{h}(e^{\frac{\alpha}{h}}-1)^{-1} e^{\alpha+2C+\frac{k\alpha}{h}-6i(\frac{\alpha}{3h}-\beta)}.$$
By the action of automorphism
 $g\mapsto
\begin{pmatrix} 0&0&1\\ 0&1&0\\ 1&0&0\end{pmatrix}
{}^{t} g^{-1}\begin{pmatrix} 0&0&1\\ 0&1&0\\ 1&0&0\end{pmatrix}$ of $G'$, which stabilizes $K$, and maps
$\begin{pmatrix} \p^{-(i+j)}&0&0\\0&\p^{-j}&0\\0&0&1\end{pmatrix}$ to $\begin{pmatrix} \p^{-(i+j)}&0&0\\0&\p^{-i}&0\\0&0&1\end{pmatrix}$, and by applying the preceding the preceding inequality to representations $\pi\circ\theta$ and $\tau\circ\theta$, we see that
$$\|c(3i,0)-c(2i,2i)\|_{V}\leq  \sqrt{C_{2}}
q^{h}(e^{\frac{\alpha}{h}}-1)^{-1} e^{\alpha+2C+\frac{k\alpha}{h}-6i(\frac{\alpha}{3h}-\beta)}.$$
where $c(0,3i)$ is invariant under the action of the subgroup $K_{1}$ of $K$ consisting of matrices of the form $\begin{pmatrix} *&*&0\\ *&*&0\\0&0&*\end{pmatrix}$,  $c(3i,0)$ invariant under $K_{2}$ of $K$ consisting of matrices of the form $\begin{pmatrix} *&0&0\\ 0&*&*\\0&*&*\end{pmatrix}$ and we have the following lemma.

   \begin{lem}

   There exists a constant $C_{3}$ (depending on $V$) such that for all $x\in V$, $y\in V$ invariant under $K_{1}$, and $z\in V$ invariant under $K_{2}$,
  we have
  $\|x\|_{V}\leq C_{3}\max(\|x-y\|_{V},\|x-z\|_{V})$.
  \end{lem}
Such subgroups $K_{1}$ and $K_{2}$ generate $K$, we have $V^{K_{1}}\cap V^{K_{2}}=0$ and the lemma follows easily. \cqfd

By applying the lemma we see that
      $$\|c(2i,2i)\|_{V}\leq C_{3} \sqrt{C_{2}}
       q^{h}(e^{\frac{\alpha}{h}}-1)^{-1} e^{\alpha+2C+\frac{k\alpha}{h}-6i(\frac{\alpha}{3h}-\beta)}.$$
It then concludes as the end of the proof of proposition \ref{prop-decr-ban-equiv'}.\cqfd

\section{Consequences}
\subsection{Extension to other groups with strong Banach property (T)}

Let $F$ be a non archimedian local field. The following theorem is theorem \ref{main} in the introduction.

\begin{thm}\label{cor-padic}
Let $G$ be a connected and almost $F$-simple algebraic group over a non archimedian local field $F$. If $\gg$ contains $sl_{3}(F)$, $G$ has strong Banach property (T).
\end{thm}

The following lemma is well-known.

\begin{lem}\label{rep-mult}
Let $\mathbb G _m$ be the multiplicative group of a field $k$. Then any algebraic representation of $G$ on a finite dimensional $k$-vector space is a finite direct sum of characters $\chi_n:\mathbb G _m\to \mathbb G _m$ defined by $x\mapsto x^n,n\in\Z$.
\end{lem}

In fact, let $\rho:\mathbb G _m\to \mathrm{GL}(V)$ be an algebraic
representation of $\mathbb G _m$ over a finite dimensional $k$-vector space $V$.
Then there exist $n_1, n_2\in\Z$ such that $\rho(x)=\sum_{n=n_1}^{n_2} P_n x^n$,
where $P_n\in \mathrm{End}(V)$. By the multiplicativity condition that
$\rho(xy)=\rho(x)\rho(y)$ we have $P_i P_j=0$ if $i\neq j$ and $P_i^2=P_i$ for
any $i,j\in\{n_1,n_1+1,..., n_2\}$. By $\rho(1)=\mathrm{Id}_V$ we have
$\sum_{n=n_1}^{n_2} P_n=Id_V$. Thus $V=\bigoplus_{n=n_1}^{n_2}V_n$, where
$V_n=\mathrm{Im}(P_n)$ and $\rho(x)v=x^n v$ for any $v\in V_n$.           
\cqfd

For an algebraic $k$-group $G$ over a field $k$, let $\textrm{Lie}(G)$ denote its Lie algebra, $\br{G}$ its points in the algebraic closure of $k$. If $G$ is semisimple and $S$ a $k$-split torus in $G$, let $\Phi(S,G)$ denote the root system associated to $(S,G)$.

The following lemma is proposition I.1.3.3 (ii) of \cite{mar}.

\begin{lem}\label{mar-uni-lie-alg}
Let $G$ be a semisimple algebraic $F$-group, $S$ a maximal $F$-split torus of $G$ contained in a maximal torus $T$. Let $\theta$ be a proper subset of the set of simple roots of $\Phi(S,G)$ with a fixed ordering. Let $V^+_\theta$(reps. $V^-_\theta$) be the $F$-subgroup of $G$, such that $\br{V^+_\theta}$(reps. $\br{V^-_\theta}$) is generated by $U_{a}$ for any root $a \in \Phi(\bar T,\bar G)$ whose restriction to $\bar S$ is a positive root in $\Phi(S,G)$ which is not a positive (reps. negative) linear combination of elements in $\theta$, where $U_a$ is the one-parameter root subgroup characterized by the existence of an isomorphism $e_a$ from the additive group $\mathbb G_a$ of the algebraic closure of $F$ to $U_a$ such that $t e_a(x)t^{-1}=e_a(a(t)x),t\in \bar T, x\in \mathbb G_a.$ Then there exist $S$-equivariant isomorphisms of $F$-varieties $\textrm{Lie}(V^+_\theta)\to V^+_\theta$ and $\textrm{Lie}(V^-_\theta)\to V^-_\theta$.
\end{lem}



The following lemma follows from proposition I.1.5.4 (III) and theorem I.2.3.1 (a) of \cite{mar}.

\begin{lem}\label{mar-gen-alg-gp}
Let $G$ be a connected, simply connected and almost $F$-simple algebraic group with $F\textrm{-rank}>0$ ($F$-isotropic), and $S,\theta,V^+_\theta,V^-_\theta$ as in lemma \ref{mar-uni-lie-alg}. Then $G$ is generated by $V^+_\theta\cup V^-_\theta$.
\end{lem}

We now begin the proof of theorem \ref{cor-padic}. By hypothesis $G$ contain a algebraic subgroup $R$ whose Lie algebra is isomorphic to $sl_{3}(F)$. Since $SL_3(F)$ is simply connected, there exists an $F$-isogeny $I:SL_3(F)\to R$, i.e. a surjective homomorphism between algebraic $F$-groups with finite kernel.

Let $\rho:F\to SL_3(F)$ be defined by $x\mapsto\begin{pmatrix}x&0&0\\0&1&0\\0&0&x^{-1}\end{pmatrix}$ for any $x\in F$, and $a=I\circ\rho(\pi)$, where $\pi$ is a uniformizer of $F$. By lemma \ref{rep-mult}, the set of eigenvalues of $Ad(a)$ is a subset of $\pi^{\Z}$ which is not $\{1\}$.
Let $S$ be a maximal $F$-split torus of $G$ containing $a$. 
We choose an ordering of $\Phi(S,G)$ such that $|b(a)|_F\leq1$ for any simple
root $b$. Let $\theta$ be the set of simple root $b$ such that $|b(a)|_F=1$, and
$V^+_\theta, V^-_\theta$ be as in lemma \ref{mar-uni-lie-alg}.        

Let $\|\cdot\|_\gg$ be the norm on $\gg$, defined by
$\|\sum_{i=1}^{\mathrm{dim}_F\gg}x_ie_i\|_\gg=\max_{1\leq i\leq
\mathrm{dim}_F\gg}|x_i|_F$, where $\{e_{i}\}_{1\leq i\leq \mathrm{dim}_F\gg}$ is
a fixed basis for $\gg$. Let $\ell'$ be a length over $G$ defined by
$\ell'(g)=\log \|Ad(g)\|_{End(\gg)}$. Let $\E$ be a class of Banach spaces
stable under duality, complex conjugation and of type $>1$. Let
$s,t,C,C'\in\R_{+}^{*}, \pro\in \CC^{\E}_{s\kappa\ell+C}(R),\pro_{m}\in
C_{c}(R)$ verify the conditions (i) and (ii) of theorem \ref{thmnonarchbanach},
where $\kappa\in\R_{+}^{*}$ such that $\ell'|_{R}\leq\kappa\ell$. Let $U$ be an
open compact subgroup of $G$ and $f=\frac{e_U}{vol(e_U)}$. Then for establishing
that $G$ has strong Banach property (T) it suffices to show that if $s$ is small
enough the series $\pro_{m}f\in C_{c}(G)$ converges in $\CC^{\E}_{s\ell'+C}(G)$
to a self adjoint idempotent $\pro'$ such that for any
$(E,\pi)\in\E_{G,s\ell'+C}$, the image of $\pi(\pro')$ consists of all
$G$-invariant vectors of $E$. First it is clear that the series $\pro_{n}f$ is a
Cauchy series in $\CC^{\E}_{s\ell'+C}(G)$ and we note $\pro'$ its limit (we see
that $\pro$ is a multiplier of $\CC^{\E}_{s\ell'+C}(G)$ and then that
$\pro'=\pro f$). Let $(E,\pi)\in \E_{G,s\ell'+C}$. It is evident that
$\pi(\pro')$ acts by identity over any $G$-invariant vector. It remains to show
that for any $x\in E$, $\pi(\pro')x$ is $G$-invariant (this will imply that
$\pro'=f^{*}\pro'=f^{*}\pro f$, so that $\pro'$ is self-adjoint). Following lemma \ref{mar-gen-alg-gp}, it suffices to show that
$\pi(\pro')x$ is $V^+_\theta$-invariant and $V^-_\theta$-invariant. Show for
example that $\pi(\pro')x$ is $V^-_\theta$-invariant. Let $V=V^-_\theta$, and
$E:\textrm{Lie}(V)\to V$ as in lemma \ref{mar-uni-lie-alg}. We know that
$\pi(\pro')x$ is fixed by $R$, then in particular by $a$. It suffices to show
that for any $Y\in\textrm{Lie}(V)$
$$\pi(E(Y))\pi(\pro')x-\pi(\pro')x=\pi(E(Y))\pi(a^{n})\pi(\pro')x-\pi(a^{n})\pi(\pro')x$$
$$=\pi(a^{n})\big(\pi(a^{-n}E(Y)a^n)-1\big)\pi(\pro')x=\pi(a^{n})\big(\pi(E(Ad(a^{-n})Y))-1\big)\pi(\pro')x$$
tends to $0$ when $n\in\N$ tends to infinity.                


Let $\textrm{Lie}(V)=\bigoplus_{\lambda\in\Lambda}\textrm{Lie}(V)^\lambda$ the decomposition of $\textrm{Lie}(V)$ under the adjoint action of $a$, $Y=\sum_{\lambda\in\Lambda} Y_\lambda$ the decomposition of $Y$, where $\Lambda\subset F$ denotes the set of eigenvalues of $Ad(a)|_{\mathrm{Lie}(V)}$. Due to the way $\theta$ is chosen, the eigenvalues of $Ad(a)|_{\textrm{Lie}(V)}$ are all of the form $\pi^{-\N^*}$. Since $U$ is an open subgroup of $G$, there exists $r>0$ such that when $Y'\in V$ and $\|Y'\|_{\textrm{Lie}(V)}\leq r$ we have $E(Y')\in U$. We put $$m=[n\kappa^{-1}\log\min_{\lambda\in\Lambda}\{|\lambda|_F\}+\kappa^{-1}\log(r/\sum_{\lambda\in\Lambda} \|Y_\lambda\|_{\textrm{Lie}(V)})],$$ where $[\cdot]$ is the integer part of a real number. When $n$ is big enough such that $m>0$, we have $$\pi(a^{n})\big(\pi(E(Ad(a^{-n})Y))-1\big)\pi(\pro_{m}f)x$$
$$= \int_{R} \pro_{m}(g)\pi(g)\big(\pi(E(Ad(g^{-1}a^{-n})Y))-1\big)\pi(f)x dg.$$
 When $$\ell'(g)\leq n\log\min_{\lambda\in\Lambda}\{|\lambda|_F\}+\log(r/\sum_{\lambda\in\Lambda}\|Y_\lambda\|_{\textrm{Lie}(V)}),$$ we have $$\|Ad(g^{-1}a^{-n})Y)\|_{\textrm{Lie}(V)}\leq e^{\ell'(g)}\sum_{\lambda\in\Lambda}|\lambda|_F^{-n}\|Y_\lambda\|_{\textrm{Lie}(V)}\leq r,$$ and hence $$\big(\pi(E(Ad(g^{-1}a^{-n})Y))-1\big)\pi(f)x=0.$$ We then have $$\pi(a^{n})\big(\pi(E(Ad(a^{-n})Y))-1\big)\pi(\pro_{m}f)x=0$$ when $n$ is big enough.

On the other hand we always have $$Ad(a^{-n})Y=\sum_{\lambda\in\Lambda} \lambda^{-n} Y_\lambda\in \bigoplus_{\lambda\in\Lambda} \O Y_\lambda,$$ where $\O$ denotes the ring of integers of $F$. Hence
$$\|\pi(a^{n})\big(\pi(E(Ad(a^{-n})Y))-1\big)\pi(\pro'-\pro_{m}f)x\|_{E}$$
 $$\leq e^{C+s\ell'(a)n} (1+C'')\|\pi(\pro'-\pro_{m}f)x\|_{E},$$
 where $C''=\sup_{t_\lambda\in \O}\|\pi(E(\sum_{\lambda\in\Lambda} t_\lambda Y_\lambda))\|_{\mathcal{L}(E)}<\infty$ depending only on $Y$. But  $$\|\pi(\pro'-\pro_{m}f)x\|_{E}\leq C' e^{-tm}
\|\pi(f)x\|_{E} $$
thanks to (ii) of theorem \ref{thmnonarchbanach} (we recall that $C'$ and $t$ are the constants of theorem \ref{thmnonarchbanach}). In total, when $n$ is big enough
$$\|\pi(a^{n})\big(\pi(E(Ad(a^{-n})Y))-1\big)\pi((\pro'-\pro_{m}f))x\|_{E}$$
 $$\leq e^{C+s\ell'(a)n} (1+C'') C' e^{-tm}
\|\pi(f)x\|_{E},$$
and if $$s<\frac{t\log \min_{\lambda\in\Lambda}\{|\lambda|_F\}}{\kappa\ell'(a)},$$ it tends to $0$ when $n$ tends to infinity.

\cqfd

We now show that strong Banach property (T) is inherited by cocompact lattices.

\begin{prop}\label{lattice}
Let $G$ be a locally compact group and $\Gamma$ a discrete cocompact subgroup of $G$. If $G$ has strong Banach property (T), so does $\Gamma$ .
\end{prop}

The existence of $\Gamma$ implies that $G$ is unimodular. Let $dg$ be the Haar
measure of $G$ such that $G/\Gamma$ has measure $1$. As $\Gamma$ is cocompact,
there exists a positive function $f\in C_{c}(G)$ such that $\sum_{\gamma\in
\Gamma}f(g\gamma)=1$ for any $g\in G$ (it implies that $\int _{G}f(g)dg=1$). Let
$X=\text{Supp}(f)$. Let $\ell$ be a length over $\Gamma$ and $\ell'$ the length
over $G$ defined by $\ell'(g)=\max\{\ell(\gamma),gX\cap X\gamma\neq
\emptyset\}$. Let $\E$ be any class of Banach spaces stable under complex
conjugation, duality and of type $>1$. Then if $\pro'\in C_{\ell'}^{\E}(G)$
verifies the conditions of definition \ref{renf-ban}, we construct an idempotent
$\pro\in C_{\ell}^{\E}(\Gamma)$ verifying the conditions of definition
\ref{renf-ban} in the following way. If $\pro'$ is the limit in
$C_{\ell'}^{\E}(G)$ of $\pro_{n}'\in C_{c}(G)$, then $\pro$ is the limit in $C_{\ell}^{\E}(\Gamma)$ of the series
$\pro_{n}\in C_{c}(\Gamma)$ defined by the following formula:
$\pro_{n}(\gamma)=\int_{G\times
G}\pro'_{n}(g_{1})f(g_{2}\gamma)f(g_{1}g_{2})dg_{1}dg_{2}$.                           

To justify it we remark that for any $(E,\pi)\in\E_{\Gamma,\ell}$, we can construct an induced representation $(E',\pi')\in\E_{G,\ell'}$ in the following way : $E'$ is the completion of the space of continuous maps $s:G\to E$ verifying $s(g\gamma)=\pi(\gamma^{-1})s(g)$ for the norm $\|s\|^{2}=\int_{G} f(g)\|s(g)\|_{E}^{2}dg$. Moreover  $\pi(\pro_{n}):E\to E$ is equal to $\beta \circ \pi'(\pro_{n}')
\circ \alpha$, where $\alpha :E\to E' $ is defined by
$\alpha(x)=(g\mapsto \sum_{\gamma\in \Gamma}f(g\gamma
)\pi(\gamma)x)$, and $\beta :E'\to E$ is given by $\beta(s)=\int_{G}f(g)s(g)dg$. The norms of $\alpha$ and $\beta$ are bounded independently of $(E,\pi)$, and the image of $\alpha$ (resp. $\beta$) of a $\Gamma$-invariant (resp. $G$-invariant) vector is $G$-invariant (resp. $\Gamma$-invariant) and if $x\in E$ is $\Gamma$-invariant, $\beta\circ\alpha(x)=x$.

Finally we calculate the involution of $\pro_n$ as defined in the introduction,
$$\pro_{n}^{*}(\gamma)=\int_{G\times G}{\pro_{n}'}^{*}(g_{1})f(g_{2}\gamma^{-1})f(g_{1}^{-1}g_{2})dg_{1}dg_{2}$$
$$=\int_{G\times G} {\pro_{n}'}^{*}(g_{1}) f(g_{3}\gamma) f(g_{3}g_{2}) dg_{1} dg_{2}$$
by the change of variables
$g_{3}=g_{1}^{-1}g_{2}\gamma^{-1}$, then as $\pro'$ is self-adjoint, so is $\pro$.

\cqfd

For non cocompact lattices, it is better to restrict to isometric representations.

\begin{prop}\label{non-ccp-latt}
Let $G$ be a locally compact group and $\Gamma$ a discrete subgroup of $G$ with finite covolume. If $G$ has Banach property (T), so does $\Gamma$.
\end{prop}

We use the notations of the proof of proposition \ref{lattice}. Let $f\in L^{1}(G)$ be positive such that $\sum_{\gamma\in \Gamma}f(g\gamma)=1$ for almost any $g\in G$ (it implies that $\int_{G}f(g)dg=1$). Then if $\pro'\in C_{0}^{\E}(G)$ verifies the conditions of definition \ref{renf-ban}, we can deduce an idempotent $\pro\in C_{0}^{\E}(\Gamma)$ verifying the conditions in definition \ref{renf-ban}. If $\pro'$ is the limit in $C_{0}^{\E}(G)$ of $\pro_{n}'\in C_{c}(G)$, then $\pro$ is the limit in $C_{0}^{\E}(\Gamma)$ of $\pro_{n}\in C_{0}^{\E}(\Gamma)$ defined by the following formula : $\pro_{n}(\gamma)=\int_{G\times G}\pro'_{n}(g_{1})f(g_{2}\gamma)f(g_{1}g_{2})dg_{1}dg_{2}$.

To justify it we remark that for any $(E,\pi)\in\E_{\Gamma,0}$, we can construct an induced representation $(E',\pi')\in\E_{G,0}$ in the following way : $E'$ is the completion of the space of continuous maps $s:G\to E$ verifying $s(g\gamma)=\pi(\gamma^{-1})s(g)$ for the norm $\|s\|^{2}=\int_{G/\Gamma} \|s(g)\|_{E}^{2}dg$. Moreover  $\pi(\pro_{n}):E\to E$ is equal to $\beta \circ \pi'(\pro_{n}')
\circ \alpha$, where $\alpha :E\to E' $ is defined by
$\alpha(x)=(g\mapsto \sum_{\gamma\in \Gamma}f(g\gamma
)\pi(\gamma)x)$, and $\beta :E'\to E$ is given by $\beta(s)=\int_{G}f(g)s(g)dg$. The norm of $\alpha$ and $\beta$ are bounded independently of $(E,\pi)$, and the image of $\alpha$ (resp. $\beta$) of a $\Gamma$-invariant (resp. $G$-invariant) vector is $G$-invariant (resp. $\Gamma$-invariant) and if $x\in E$ is $\Gamma$-invariant, $\beta\circ\alpha(x)=x$.

The same calculation at the end of the demonstration of proposition \ref{lattice} shows that $\pro$ is self-adjoint.

\cqfd

\subsection{Embedding of expanders}

Let $F$ be a non archimedian local field. Let $G$ be a connected almost $F$-simple algebraic group over $F$ whose Lie algebra contains $sl_{3}(F)$. Let $\Gamma$ be a lattice of $G(F)$. Following corollary \ref{cor-padic} and proposition \ref{non-ccp-latt}, $\Gamma$ has Banach property (T).

Let $(\Gamma_{i})_{i\in\N}$ be a series of subgroups of $\Gamma$ such that $\sharp(\Gamma/\Gamma_{i})$ tends to infinity (we know that such a series exits). We choose a finite symmetric system of generators $S$ of $\Gamma$. We suppose $S^{2}\cap\Gamma_{i}=\{1\}$ for any $i$. The system $S$ equips $X_{i}=\Gamma/\Gamma_{i}$ with a graph structure and we note $d_{i}$ its associated metric. As $\Gamma$ has the usual property (T), $X_{i}$ forms a family of expanders.

We say that the series $X_{i}$ is embedded uniformly in a Banach space $E$, if there exists a function $\rho:\N\to\R_{+}$ tends to $+\infty$ at infinity and $1$-Lipschitz maps $f_{i}:X_{i}\to E$ such that $\|f_{i}(x)-f_{i}(y)\|_{E}\geq\rho(d_{i}(x,y))$ for any $i\in\N$ and $x,y\in X_{i}$.

\begin{thm}\label{expanders}
The series of expanders $(X_{i},d_{i})$ does not admit a uniform embedding in any Banach space of type $>1$.
\end{thm}

In fact any subseries does not admit such a uniform embedding. We have the following more precise result.

\begin{prop}\label{poincare-ineq}
Let $E$ be a Banach space of type $>1$. There exists a constant $C$ such that for any $i\in\N$, and for any map $f:X_{i}\to E$, we have
$$\esp{x,y\in X_{i}}\|f(x)-f(y)\|^{2}\leq C \esp{x,y\text{ neighbors in } X_{i}}\|f(x)-f(y)\|^{2}.$$
\end{prop}

It is well known that proposition \ref{poincare-ineq} implies the theorem. If
the series $(X_{i},d_{i})$ admitted a uniform embedding in a Banach space $E$ of
type $>1$, we would have $1$-Lipschitz maps $f_{i}:X_{i}\to E$ and a map
$\rho:\N\to\R_{+}$ tending to $+\infty$ such that for any $i\in\N$ and $x,y\in
X_{i},$ we have $$d_i(x,y)\geq\|f_{i}(x)-f_{i}(y)\|\geq \rho(d_i(x,y)).$$
 
 We have $\esp{x,y\in
X_{i}}\rho(d_i(x,y))^{2}$ tends to $+\infty$ when $i$ tends to infinity. In
fact, for any $k\in\N$ we note $K_i(k)$ the set of elements $(x,y)\in (X_i)^2$ such
that $d_i(x,y)\leq k,$ then we have $ K_i(k)\leq (\sharp S)^k \sharp X_i.$
For any big enough $N\in \N,$ let $k_N\in\N$ such that $a\geq k_N$ implies
$\rho(a)\geq N.$ Then we have$$\esp{x,y\in X_{i}}\rho(d(x,y))^{2}\geq
\frac{1}{(\sharp X_i)^2}M^2((\sharp X_i)^2-K_i(k_N))\geq N^2(1-\frac{(\sharp
S)^{k_N}}{\sharp X_i}).$$ We see the right handside of the above inequality
tends to infinity when $i$ tends to infinity.

However, this contradicts with the inequality in proposition \ref{poincare-ineq}
since $$\esp{x,y\text{ neighbors in } X_{i}}\|f(x)-f(y)\|^{2}\leq \esp{x,y\text{
neighbors in } X_{i}}d_i(x,y)^{2}=1.$$
\cqfd

It remains to show proposition \ref{poincare-ineq} We note $E_{i}$ the space of functions of $X_{i}$ into $E$, endowed with the norm $\|f\|_{E_{i}}^{2}=\esp{x\in X_{i}}\|f(x)\|^{2}_{E}$. There exits a class of Banach spaces $\E$ of type $>1$ that contains $E_{i}$. In fact any class of type $>1$ is included in a class of type $>1$ stable under direct $\ell^{2}$ sum. As $E_{i}$ is an isometric representation of $\Gamma$ we have $E_{i}\in\E_{\Gamma,0}$. We recall that $C_{0}^{\E}(\Gamma)$ is the Banach algebra of $C_{c}(\Gamma)$ under the completion with norm $\|f\|=\sup_{(E,\pi)\in\E_{\Gamma,0}}\|\pi(f)\|_{\L(E)}$. As $\Gamma$ has Banach property (T), there exits an idempotent $\pro \in C_{0}^{\E}(\Gamma)$ such that for any representation $(E,\pi)$ in the class $\E_{\Gamma,0}$, the image of $\pi(\pro)$ consists of exactly $\Gamma$-invariant vectors in $E$.

The $\Gamma$-invariant vectors in $E_{i}$ are exactly the constant functions of $X_{i}$ in $E$. For any function $f\in E_{i}$, we have $\pro f=m_{f}$, where $m_{f}=\esp{x\in X_{i}}f(x)\in E$ is the average of $f$, considered as a constant function over $X_{i}$. It is well known that
$$\esp{x,y\in X_{i}}\|f(x)-f(y)\|^{2}=2\esp{x\in X_{i}}\|f(x)-m_{f}\|^{2}=2\|f-\pro f\|_{E_{i}}^{2}.$$

Let $\pro_{1}\in C_{c}(\Gamma)$, of integral $1$, such that $\|\pro-\pro_{1}\|_{\CC^{\E}_{0}(\Gamma)}\leq \frac{1}{2}$. Then $$\|\pro f-\pro_{1}f\|_{E_{i}}=\|(\pro-\pro_{1})(f-m_{f})\|_{E_{i}}\leq \frac{1}{2}\|f-m_{f}\|_{E_{i}}= \frac{1}{2}\|f-\pro f\|_{E_{i}}$$
and $\|f-\pro f\|_{E_{i}}\leq 2\|f-\pro_{1}f\|_{E_{i}}$. 

We see
that there exits a constant $C_{1}$ (depending only on $\pro_{1}$) such that
$$\|f-\pro_{1}f\|_{E_{i}}^{2}\leq  C_{1}\esp{x,y\text{ neighbors in }
X_{i}}\|f(x)-f(y)\|^{2}.$$ In fact, For any $\gamma\in S^k,$ by triangular
inequality we have $$\|f-\gamma f\|_{E_i}\leq k\big(\sum_{s\in
S}\|f-sf\|_{E_i}\big).$$ If $\mathrm{supp}(\pro_1)\in S^k,$ we see that
$$\|f-\pro_1f\|_{E_i}\leq\sum_{\gamma\in\Gamma}|\pro_1(\gamma)|k\sum_{s\in
S}\|f-s f\|_{E_i}$$$$=\Big(\sum_{\gamma\in\Gamma}|\pro_1(\gamma)|k\sharp
S\Big)\esp{x,y\text{ neighbors in } X_{i}}\|f(x)-f(y)\|^{2}.$$
\cqfd                  

\subsection{Fixed point property for affine actions}

Corollary \ref{prop-f} of the introduction is a consequence of the following proposition and of corollary \ref{cor-padic}.

\begin{prop}\label{fixed-point}
Let $G$ be a locally compact group. If $G$ has strong Banach property (T) in the sense of definition \ref{renf-ban}, then $G$ has Banach property (F) in the sense of definition \ref{def-f}.
\end{prop}

In fact let $E$ be an affine Banach space whose underlying vector Banach space $E_{0}$ is of type $>1$ and let $\rho$ be a continuous action of $G$ on $E$ by affine isometry. Let $x_{0}\in E$ be a point and $\ell:G\to \R_{+}$ a length defined by $\ell(g)=\|\rho(g)(x_{0})-x_{0}\|$. We put $\tilde{E}=E_{0}\oplus\C$ (with the norm of the $\ell^{2}$ direct sum) and we identify $E$ with the hyperplane $E_{0}\times\{1\}$ of $\tilde{E}$ such that $x_{0}$ is sent to $(0,1)$. Let $\pi$ be the linear representation of $G$ on $\tilde{E}$ such that $G$ preserves the hyperplane $E_{0}\times\{1\}$ and acts on it by $\rho$, through the previous identification. The representation $\pi$ is not isometric but we have $\|\pi(g)\|\leq 1+\ell(g)$. For any $s>0$ there exists $C$ such that $1+\ell(g)\leq Ce^{s\ell(g)}$ for any $g\in G$. As $G$ has strong Banach property (T), and we have a $G$-equivariant surjection $\tilde{E}\to\C$ (where $\C$ is endowed with the trivial representation), by the argument in the second remark after definition \ref{renf-ban}, we see that there exits a $G$-invariant vector in $\tilde{E}$ whose image in $\C$ is equal to $1$. This means exactly that $\rho$ has a fixed point.

{\bf Remark} More generally, if $G$ has strong Banach property (T), for any class $\E$ of type $>1$ and for any length $\ell$ on $G$ there exits $s>0$ such that for any $C\in\R_{+}$, any affine action of $G$ over an affine Banach space whose underlying vector Banach space belongs to $\E_{G,C+s\ell}$ admits a fixed point.

\end{document}